\documentclass{amsart}

\usepackage{amssymb}
\usepackage{txfonts}
\usepackage{cleveref}
\usepackage{microtype}
\usepackage{mathrsfs}
\usepackage{graphicx}
\usepackage{cite}
\usepackage{xcolor}

\newtheorem{theorem}{Theorem}[section]
\newtheorem{lemma}[theorem]{Lemma}
\newtheorem{corollary}[theorem]{Corollary}

\theoremstyle{definition}
\newtheorem{definition}[theorem]{Definition}
\newtheorem{example}[theorem]{Example}
\newtheorem{proposition}[theorem]{Proposition}

\theoremstyle{remark}
\newtheorem{remark}[theorem]{Remark}

\numberwithin{equation}{section}

\begin{document}

\title{Subshifts of finite type and matching for intermediate $\beta$-transformations}


\author[Y.\ Sun]{Yun Sun}
\address[Yun Sun]{School of Mathematic, South China University of Technology, Guangzhou, 510461, China}
\email{masy2021@mail.scut.edu.cn}
\author[B.\ Li]{Bing Li}
\address[Bing Li]{School of Mathematic, South China University of Technology, Guangzhou, 510461, China}
\email{scbingli@scut.edu.cn}
\author[Y M.\ Ding]{Yiming Ding$^{\dag}$}

\address[Yiming Ding]{College of Science, Wuhan University of Science and Technology, Wuhan, 430081, China}
\email{dingym@wust.edu.cn}



\thanks{ $^{\dag}$Corresponding author. Y M. Ding was partially supported by NSFC 12271418. B. Li was partially supported by NSFC 12271176. \\
}

\subjclass[2010]{Primary: 37E05, 37B10; Secondary: 11A67, 11R06}

\date{\today}

\dedicatory{}


\begin{abstract}
We focus on the relationships between matching and subshift of finite type for intermediate $\beta$-transformations $T_{\beta,\alpha}(x)=\beta x+\alpha $ ($\bmod$ 1), where $x\in[0,1]$ and $(\beta,\alpha) \in \Delta:= \{ (\beta, \alpha) \in \mathbb{R}^{2}:\beta \in (1, 2) \; \rm{and} \; 0 < \alpha <2 - \beta\}$. We prove that if the kneading space $\Omega_{\beta,\alpha}$ is a subshift of finite type, then $T_{\beta,\alpha}$ has matching. Moreover, each $(\beta,\alpha)\in\Delta$ with $T_{\beta,\alpha}$ has matching corresponds to a matching interval, and there are at most countable different matching intervals on the fiber. Using combinatorial approach, we construct a pair of linearizable periodic kneading invariants
and show that, for any $\epsilon>0$ and $(\beta,\alpha)\in\Delta$ with $T_{\beta,\alpha}$ has matching, there exists $(\beta,\alpha^{\prime})$ on the fiber with $|\alpha-\alpha^{\prime}|<\epsilon$, such that $\Omega_{\beta,\alpha^{\prime}}$ is a subshift of finite type. As a result, the set of $(\beta,\alpha)$ for which $\Omega_{\beta,\alpha}$ is a subshift of finite type is dense on the fiber if and only if the set of $(\beta,\alpha)$ for which $T_{\beta,\alpha}$ has matching is dense on the fiber.


\par Key words: $\beta$-transformation; subshifts of finite type; matching; kneading invariants;
\end{abstract}

\maketitle


\section{Introduction and statement of main results}
\subsection{Introduction}
\
\par Since the pioneering work of R\'enyi \cite{renyi1957representations} and Parry \cite{parry1960beta,parry1964representations}, intermediate $\beta$-shifts have been well-studied by many authors, and have been shown to have important connections with ergodic theory, fractal geometry and number theory.  We refer the reader to \cite{blanchard1989,komornik2011,sidorov2003} for articles concerning this topic.

For $\beta > 1$ and $x \in [0,1/(\beta-1)]$, a word $(\omega_{n})_{n \in \mathbb{N}}$ with symbols in the alphabet $\{ 0, 1\}$ is called a $\beta$-\textsl{expansion} of $x$ if
$x = \sum_{k= 1}^{\infty} \omega_{k} \ \beta^{-k}$. Through iterating the maps $G_{\beta} \colon x \mapsto \beta x \bmod 1$ on $[0,1)$ and $L_{\beta} \colon x \mapsto \beta (x - 1) \bmod 1$ on $(0,1]$ with $\beta \in (1,2]$, see Figure 1, one obtains subsets of $\{0, 1\}^{\mathbb{N}}$ known as the greedy and (normalised) lazy $\beta$-shifts, respectively. Each point $\omega^{+}$ of the greedy $\beta$-shift is a $\beta$-expansion, and corresponds to a unique point in $[0, 1]$. Similarly, each point $\omega^{-}$ of the lazy $\beta$-shift is a $\beta$-expansion, and corresponds to a unique point in $[(2-\beta)/(\beta - 1), 1/(\beta - 1)]$.  Note that, if $\omega^{+}$ and $\omega^{-}$ are $\beta$-expansions of the same point, then $\omega^{+}$ and $\omega^{-}$ do not need to be equal, see for instance \cite[Theorem 1]{komornik1998}.
\par Denote $\Delta := \{ (\beta, \alpha) \in \mathbb{R}^{2} \colon \beta \in (1, 2) \; \text{and} \; 0 <\alpha <2 - \beta\}$. The intermediate $\beta$-shifts $\Omega_{\beta, \alpha}$ arises from the intermediate $\beta$-transformations $T^{\pm}_{\beta,\alpha} \colon [0,1] \circlearrowleft$, where $(\beta,\alpha)\in\Delta$ and the maps $T_{\beta, \alpha}^{\pm}$ are defined as follows.  Let $c = c _{\beta, \alpha} := (1-\alpha)/\beta$, set
\begin{align*}
T^{+}_{\beta, \alpha}(x) :=
\begin{cases}
\beta x + \alpha \ (\bmod 1) & \text{if} \; x \neq c,\\
0 & \text{if} \; x = c,
\end{cases}
\quad \text{and} \quad
T^{-}_{\beta, \alpha}(x) :=
\begin{cases}
\beta x + \alpha \ (\bmod 1) & \text{if} \; x \neq c,\\
1 & \text{if} \; x = c.
\end{cases}
\end{align*}
The maps $T_{\beta, \alpha}^{\pm}$ are equal everywhere except at the point $c$ and $T^{-}_{\beta, \alpha} (x) = 1 - T_{\beta, 2 - \beta - \alpha}^+(1 - x)$, for all $x \in [0,1]$, see Figure 2. Notice, when $\alpha = 0$, the maps $G_{\beta}$ and $T^{+}_{\beta, \alpha}$ coincide on $[0,1)$, and when $\alpha = 2-\beta$, the maps $L_{\beta}$ and $T^{-}_{\beta, \alpha}$ coincide on $(0,1]$. Observe that, for all $(\beta,\alpha)\in\Delta$, the symbolic space $\Omega_{\beta,\alpha}$ of $T^{+}_{\beta,\alpha}$ and $T^{-}_{\beta,\alpha}$, see Section 2.2 for a formal definition, is always a subshift, meaning that it is invariant under the left shift map. Subshifts of finite type play an essential role in the study of dynamical systems, and can be completely described by a finite set of forbidden words (see Section 2.1). For convenience, we write subshift of finite type as \textsl{SFT}. If $\Omega\subseteq\{0,1\}^{\mathbb{N}}$ is a factor of a SFT, then it is called \textsl{sofic}. SFT has a simple representation as a finite directed graph, thus dynamical and combinatorial questions about SFT can be phrased in terms of an adjacency matrix, making them more tractable. Hence it is of interest to classify the points in $\Delta$ for which $\Omega_{\beta,\alpha}$ is a SFT.

\par Given $(\beta, \alpha) \in \Delta$, the $\beta$-expansions of the critical point $c$ given by $T_{\beta, \alpha}^{\pm}$ are called the kneading invariants of $\Omega_{\beta, \alpha}$ (see Section \ref{susbsection2.2}), denoted as $(k_{+},k_{-})$. It was mentioned in \cite[Theorem 2]{hubbard1990classification} that the kneading invariants completely determine $\Omega_{\beta, \alpha}$. The previous results for the greedy, lazy and intermediate $\beta$-shifts to be SFT refer to Theorem \ref{lazysft} and Theorem \ref{intersft} in Section \ref{susbsection2.2}, these results immediately give us that the set of parameters in $\Delta$ which give rise to SFT or sofic is countable. In another article \cite{li2019denseness} by Li et al., it was proved that the set of $(\beta,\alpha)$ belonging to $\Delta$ for which $\Omega_{\beta,\alpha}$ is a SFT is dense in $\Delta$. Denote the fiber $$\Delta(\beta):=\{(\beta,\alpha)\in\mathbb{R}^{2}:0<\alpha<2-\beta\},$$
where $\beta\in(1,2)$ is fixed. A natural question arises that, whether the set of $(\beta,\alpha)$
in $\Delta(\beta)$ with $\Omega_{\beta,\alpha}$ being a SFT is dense in $\Delta(\beta)$? The article \cite{quackenbush2020periodic} gives a positive answer to this question in the case that $\beta$ is a multinacci number, which is the unique real solution to the
equation $\beta^{k}=\beta^{k-1}+\cdots+\beta+1$ $(2\leq k<\infty)$ in interval $(1, \ 2)$. How about other Pisot numbers? Here we give an equal characterization in Theorem \ref{th1} via matching.

\par We say intermediate $\beta$-transformation $T_{\beta,\alpha}$ has matching, if there exist a finite integer $n$ such that $(T_{\beta,\alpha}^{+})^{n}(0)=(T_{\beta,\alpha}^{-})^{n}(1)$, and the smallest such $n$ is referred to be matching time. Matching has attracted attentions in the study of iterated piecewise maps and is often related with Markov partitions, entropy and invariant measures. Combining with the result in \cite{parry1960beta} about $T_{\beta,\alpha}$-invariant density, matching then immediately implies the $T_{\beta,\alpha}$-invariant density is piecewise constant. Note that we do not consider those $\beta$-transformations with fixed points, that is, $\alpha=0$ or $2-\beta$. For instance, let $\beta$ be the golden mean, then both $\Omega_{\beta,0}$ and $\Omega_{\beta,2-\beta}$ are subshifts of finite type, and their $T_{\beta,\alpha}$-invariant densities are piecewise constant.
But by the definition of matching, $T_{\beta,0}$ and $T_{\beta,2-\beta}$ do not have matching. In several parametrised families where matching was studied, it turned out that matching occurs prevalently. For example, it is shown in \cite{kraaikamp2012natural} that the set of $\alpha$'s for which the $\alpha$-continued fraction map has matching, has full Lebesgue measure. However, for piecewise linear transformation, prevalent matching appears to be rare. Such as generalised $\beta$-transforamtion (that is $\beta>1$ and $\alpha\in[0,1]$), it was proved in \cite{bruin2017} that when $\beta$ is a quadratic Pisot number or a tribonacci number, matching occurs on subset of $\alpha\in[0,1]$ with full Lebesgue measure, and the authors conjecture that the result holds for $\beta$ being a Pisot number.
\par Here we focus on the intermediate $\beta$-transformations, which have only two branches. We say $T_{\beta,\alpha}$ and $T_{\beta,\alpha^{\prime}}$ have the same matching if matching occurs at the same time $n$, and the first $n+1$ symbols of their kneading invariants are identical. For example, let $(k_{+},k_{-})=(100(10)^{\infty},011(10)^{\infty})$ and $(k^{\prime}_{+},k^{\prime}_{-})=(100(010)^{\infty},011(010)^{\infty})$, then the intermediate $\beta$-transformations correspond to $(k_{+},k_{-})$ and $(k^{\prime}_{+},k^{\prime}_{-})$ have the same matching. It was proved in  \cite[Proposition 5.1]{bruin2017} that matching occurs on the whole fiber $\Delta(\beta)$ when $\beta$ is a multinacci number. This gives that the set of parameters in $\Delta$ which give rise to matching is uncountable. Moreover, multinacci number is so special such that all the $(\beta,\alpha)\in\Delta(\beta)$ have the same matching, and 
we prove that, for $\beta$ being any other Pisot number which is not multinacci number, the whole fiber will not have the same matching.
This will not contradict with the conjecture in \cite{bruin2017} since there are at most countable different matching intervals (see the proof of Proposition \ref{SFTIFF}) and their collection may be of full Lebesgue measure. Notice that both the result in \cite{quackenbush2020periodic} about SFT and the result in \cite{bruin2017} about matching are related with multinacci number, and obtain good results.  Natural questions are that, what is the relationship between two different properties? what if consider other Pisot numbers? In this paper, we mainly focus on the relationships between matching and SFT, see the following main results.

\subsection{Statement of main results}
\
\par Let $A\subset \Delta$, denote $\overline{A}$ as the closure under the usual Euclidean metric. Before stating main results, we list some useful notations.
\begin{align*}
\begin{cases}
\mathcal{M}:=\{(\beta,\alpha)\in\Delta: T_{\beta,\alpha} \ {\rm has  \ matching} \},\\
\mathcal{F}:=\{(\beta,\alpha)\in\Delta:\Omega_{\beta,\alpha} \ {\rm is  \ a \ SFT} \},\\
\mathcal{S}:=\{(\beta,\alpha)\in\Delta:\Omega_{\beta,\alpha} \ {\rm is   \ sofic} \},\\
\mathcal{M}(\beta):=\Delta(\beta)\cap\mathcal{M},\\
\mathcal{F}(\beta):=\Delta(\beta)\cap\mathcal{F},\\
\mathcal{S}(\beta):=\Delta(\beta)\cap\mathcal{S},\\
I(\beta,\alpha):= \{(\beta,\alpha^{\prime})\in\Delta(\beta): T_{\beta,\alpha^{\prime}} \ {\rm and} \  T_{\beta,\alpha}\ {\rm have\ the \ same \ matching}\},\\
\mathcal{F}(\beta,\alpha):= I(\beta,\alpha)\cap\mathcal{F}.
\end{cases}
\end{align*}

\begin{theorem}\label{th1}
\
\par
\begin{enumerate}
 \item $\mathcal{F}\subsetneq \mathcal{M}$.
\item $\overline{\mathcal{F}(\beta)}=\overline{\mathcal{M}(\beta)}$.
\end{enumerate}
\end{theorem}
As a result, $\mathcal{F}(\beta)$ is dense in $\Delta(\beta)$ if and only if $\mathcal{M}(\beta)$ is dense in $\Delta(\beta)$. Let $(\beta,\alpha)\in\mathcal{M}$, the sufficient and necessary condition for $(\beta,\alpha)\in\mathcal{F}$ is stated in Proposition \ref{SFTIFF}. Furthermore, there exist examples indicate that neithor $\mathcal{M}\setminus \mathcal{S}$ nor $\mathcal{S}\setminus \mathcal{M}$ is empty.


\begin{remark}\label{emptyset}
\
\par
\begin{enumerate}
 \item $\overline{\mathcal{M}}=\overline{\mathcal{F}}=\overline{\mathcal{S}}=\overline{\Delta}$.
\item $\mathcal{F}(\beta)=\emptyset$ if and only if $\mathcal{M}(\beta)=\emptyset$.
\item Even $\beta$ is a Perron number, $\mathcal{M}(\beta)$ may be empty, see Example \ref{emptybeta}.
\end{enumerate}
\end{remark}

\begin{theorem}\label{th2}
Let $(\beta,\alpha)\in\mathcal{M}$ and $(T_{\beta,\alpha}^{+})^{m-1}(0)=(T_{\beta,\alpha}^{-})^{m-1}(1)$. Write $k_{+}=(10a_{3}\cdots a_{m}\cdots)$ and $k_{-}=(01b_{3}\cdots b_{m}\cdots)$. Then
\begin{enumerate}
\item $I(\beta,\alpha)$ is a subinterval of $\Delta(\beta)$.
\item $\overline{\mathcal{F}(\beta,\alpha)}=\overline{I(\beta,\alpha)}$.
\item Denote $(\beta,\alpha_{l})$ and $(\beta,\alpha_{r})$ as the left and right endpoints of $I(\beta,\alpha)$,
\begin{enumerate}
\item if $a_{m}=0$ and $b_{m}=1$, then $\{(\beta,\alpha_{l}),(\beta,\alpha_{r})\}\subset\mathcal{F}(\beta)$;
\item if $a_{m}=1$ and $b_{m}=0$, then $\{(\beta,\alpha_{l}),(\beta,\alpha_{r})\}\subset\mathcal{S}(\beta)\setminus\mathcal{M}(\beta)$.
\end{enumerate}
\end{enumerate}
\end{theorem}
We also call $I(\beta,\alpha)$ a matching interval. Let $(\beta,\alpha)\in\mathcal{M}$, we give two ways to calculate $I(\beta,\alpha)$, see the proof of Theorem \ref{th2} (1) and Lemma \ref{calendpoint}. Moreover, if consider the closure of $I(\beta,\alpha)$, whether the endpoints of $I(\beta,\alpha)$ belong to $\mathcal{F}$ relies on the value of $a_{m}$ and $b_{m}$.

\begin{remark}\label{tails}
\
\par
\begin{enumerate}
 \item When $I(\beta,\alpha)$ is a singleton,  $\mathcal{F}(\beta,\alpha)=I(\beta,\alpha)=\{(\beta,\alpha)\}$.
 \item There are at most countable different matching intervals on $\Delta(\beta)$.
 \item Except for the cases $(\beta,\alpha_{l})=(\beta,0)$ and $(\beta,\alpha_{r})=(\beta,2-\beta)$:
\begin{enumerate}
\item $k_{-}$ of $\Omega_{\beta,\alpha_{l}}$ is periodic, $k_{+}$ of $\Omega_{\beta,\alpha_{r}}$ is periodic.
\item if $a_{m}=0$ and $b_{m}=1$, $I(\beta,\alpha)$ is closed; if $a_{m}=1$ and $b_{m}=0$, $I(\beta,\alpha)$ is open.
\end{enumerate}
\end{enumerate}

\end{remark}

\begin{corollary}\label{cor1.4} For any $I(\beta,\alpha_{1})\neq I(\beta,\alpha_{2})$, we have $\overline{I(\beta,\alpha_{1})}\cap \overline{I(\beta,\alpha_{2})}=\emptyset$.

\end{corollary}

As a result, the right endpoint of a matching interval can not be the left endpoint of another matching interval on the fiber $\Delta(\beta)$.
The conjecture in \cite{bruin2017} stated as, $\mathcal{M}(\beta)$ is of full Lebesgue measure if $\beta$ is a Pisot number. Corollary \ref{cor1.4} does not contradict with this conjecture since there may still have a lot of $(\beta,\alpha)$ belonging to $\mathcal{M}$ between any two matching intervals.
\begin{corollary}\label{cor1.3} Let $\alpha\in(0,2-\beta)$. $I(\beta,\alpha)=\Delta(\beta)$ if and only if $\beta$ is a multinacci number.

\end{corollary}
It was proved in \cite[Proposition 5.1]{bruin2017} and \cite[Proposition 1]{quackenbush2020periodic} that when $\beta$ is a multinacci number, $I(\beta,\alpha)=\Delta(\beta)$ for any $\alpha\in(0,2-\beta)$. Here we show that if $\beta$ is not a multinacci number, then for any $\alpha\in(0,2-\beta)$,  $I(\beta,\alpha)\neq\Delta(\beta)$.

\section{Preliminaries}

\subsection{Subshifts}
\
\par We equip the space $\{0,1\}^\mathbb{N}$ of infinite sequences with the topology induced by the usual metric $d \colon \{0,1\}^\mathbb{N} \times \{0,1\}^\mathbb{N} \to \mathbb{R}$ which is given by
\begin{align*}
d(\omega, \nu) \coloneqq
\begin{cases}
0 & \text{if} \; \omega = \nu,\\
2^{- \lvert\omega \wedge \nu\rvert + 1} & \text{otherwise}.
\end{cases}
\end{align*}
Here $\rvert \omega \wedge \nu \lvert \coloneqq \min \, \{ \, n \in \mathbb{N} \colon \omega_{n} \neq \nu_n \}$, for all $\omega = (\omega_{1}\omega_{2}\dots) , \nu = ( \nu_{1} \nu_{2}\dots) \in \{0, 1\}^{\mathbb{N}}$ with $\omega \neq \nu$. Note that the topology induced by $d$ on $\{ 0, 1\}^{\mathbb{N}}$ coincides with the product topology on $\{ 0, 1\}^{\mathbb{N}}$.  For $n\in\mathbb{N}$ and $\omega\in\{0,1\}^{\mathbb{N}}$, we set $\omega|_{1}^{n}=\omega|_{n}=(w_{1}\cdots w_{n})$ and call $n$ the length of $\omega|_{n}$. We let $\sigma \colon \{ 0, 1 \}^{\mathbb{N}} \circlearrowleft$ denote the \textsl{left-shift map} which is defined by $\sigma(\omega_{1} \omega_{2} \dots) \coloneqq (\omega_{2} \omega_{3}\dots)$.  A \textsl{subshift} is any closed subset $\Omega \subseteq \{0,1\}^\mathbb{N}$ such that $\sigma(\Omega) \subseteq \Omega$.  Given a subshift $\Omega$ and $n \in \mathbb{N}$ we set
\begin{align*}
\Omega\lvert_{n} \coloneqq \left\{ (\omega_{1} \dots \omega_{n}) \in \{ 0, 1\}^{n} \colon \text{there exists} \, \omega \in \Omega \, \ \text{with} \ \, \omega|_{n} = (\omega_{1} \dots \omega_{n}) \right\}
\end{align*}
and denote by $\Omega^{*} \coloneqq \bigcup_{n = 1}^{\infty} \Omega\lvert_{n}$ for the collection of all finite words. For $\xi\in\Omega^{*}$, we denote $|\xi|$ as the length of $\xi$. And we denote by $\#\Omega|_{n}$ the cardinality of $\Omega|_{n}$. A subshift $\Omega$ is said to be \textsl{of finite type} if there exists a finite set $F$ of finite words such that
\par \  (i) $\nu|_{n}\notin F$ for all $\nu\in\Omega$ and $n\in\mathbb{N}$;
\par \  (ii) if $\nu\in \{0,1\}^{\mathbb{N}}\setminus\Omega$, then there exist integers $n>0$ and $m\geq0$ such that $\sigma^{m}(\nu)|_{n}\in F$.
For convenience, we write it as \textsl{SFT}. The set $F$ is often referred to as the set of \textsl{forbidden words} of $\Omega$.
\par For $n, m \in \mathbb{N}$ and $\nu = (\nu_{1}\dots \nu_{n}),\, \xi = (\xi_{1} \dots \xi_{m}) \in \{ 0, 1\}^{*}$, set
$
\nu \xi\coloneqq (\nu_{1} \dots \nu_{n} \xi_{1} \dots\xi_{m});
$
we use the same notation when $\xi \in \{ 0, 1\}^{\mathbb{N}}$.  An infinite word $\omega = (\omega_{1} \omega_{2} \dots) \in \{0, 1\}^{\mathbb{N}}$ is called \textsl{periodic} with \textsl{period} $n \in \mathbb{N}$ if and only if, $(\omega_{1}, \dots, \omega_{n}) = (\omega_{(m - 1)n + 1} \dots \omega_{m n})$, for all $m \in \mathbb{N}$; in which case we write $\omega = (\omega_{1} \dots \omega_{n})^{\infty}$.  Similarly,  $\omega = (\omega_{1}\omega_{2} \dots) \in \{0, 1\}^{\mathbb{N}}$ is called \textsl{eventually periodic} with \textsl{period} $n \in \mathbb{N}$ if and only if there exists $k \in \mathbb{N}$ such that $(\omega_{k+1} \dots \omega_{k+n}) = (\omega_{k+(m - 1)n + 1} \dots, \omega_{k+ m n})$ for all $m \in \mathbb{N}$; in which case we write $\omega = \omega_{1}\dots \omega_{k} (\omega_{k+1} \dots \omega_{k+n})^{\infty}$.

\begin{figure}[t]
\begin{center}
\includegraphics[width=0.35\textwidth]{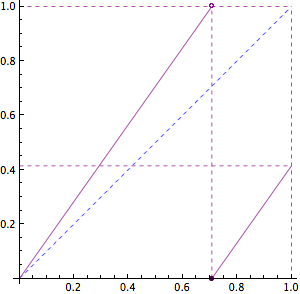}
\hspace{4em}
\includegraphics[width=0.35\textwidth]{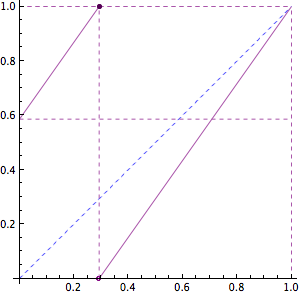}
\caption{$\beta = 2^{1/2}$. Left: $G_{\beta} = T^{+}_{\beta, 0}$. Right: $L_{\beta} = T^{-}_{\beta, 2 -  \beta}$.$G_{\beta} = T^{+}_{\beta, 0}$.}
\label{fig:fig1}
\vspace{1em}
\includegraphics[width=0.35\textwidth]{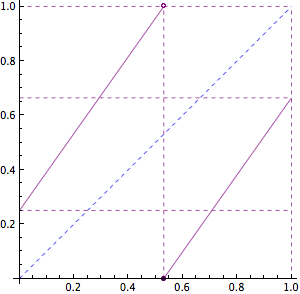}
\hspace{4em}
\includegraphics[width=0.35\textwidth]{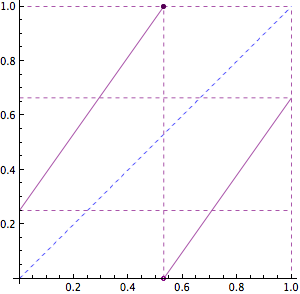}
\caption{ $\beta = 2^{1/2}$ and $\alpha = 1/4$. Left: $T^{+}_{\beta, \alpha}$. Right: $T^{-}_{\beta, \alpha}$.}
\label{fig:fig2}
\end{center}
\end{figure}

\subsection{Intermediate $\beta$-shifts of finite type}\label{susbsection2.2}
\
\par Throughout this section let $(\beta, \alpha) \in \Delta$ be fixed and the critical point $c =(1 - \alpha)/\beta$. The \textsl{$T^{\pm}_{\beta, \alpha}$-expansion} $\tau_{\beta, \alpha}^{\pm}(x)$ of $x \in [0, 1]$ is defined to be the word $(\omega^{\pm}_{1}\omega^{\pm}_{2} \dots ) \in \{0 ,1\}^{\mathbb{N}}$, where, for $n \in \mathbb{N}$,
\begin{align*}
\omega^{+}_{n} \coloneqq \begin{cases}
0 & \quad \text{if } (T^{+}_{\beta,\alpha})^{n-1}(x) < c,\\
1 & \quad \text{otherwise,}
\end{cases}
\quad \text{and} \quad
\omega^{-}_{n} \coloneqq \begin{cases}
0 & \quad \text{if } \, (T^{-}_{\beta,\alpha})^{n-1}(x) \leq c,\\
1 & \quad \text{otherwise.}
\end{cases}
\end{align*}
We will denote the images of the unit interval under $\tau_{\beta, \alpha}^{\pm}$ by $\Omega^{\pm}_{\beta, \alpha}$, respectively, and set $\Omega_{\beta, \alpha} \coloneqq \Omega_{\beta, \alpha}^{+} \cup \Omega_{\beta, \alpha}^{-}$.  The \textsl{kneading invariants} of $\Omega_{\beta,\alpha}$ is defined to be the pair of sequences $(k_{+},k_{-})=(\tau^{+}_{\beta, \alpha}(c),\tau^{-}_{\beta, \alpha}(c))$.

\begin{remark}\label{rem:00..11..}
 We have $k_{+} = 10^{\infty}$ if and only if $\alpha = 0$; and $k_{-} =01^{\infty}$ if and only if $\alpha = 2-\beta$.  Moreover, $k(0)=\tau_{\beta, \alpha}^{+}(0) = \sigma(k_{+})$ and $k(1)=\tau_{\beta, \alpha}^{-}(1) = \sigma(k_{-})$.
\end{remark}
\begin{theorem}[{\cite[Proposition 4]{alseda1996},\cite[Theorem 5.1]{barnsley2012} and \cite[Theorem 2]{kalle2012beta}}]\label{thm:Structure}
For $(\beta, \alpha) \in \Delta$, the kneading spaces $\Omega_{\beta, \alpha}^{\pm}$ are completely determined by the kneading invariants of $\Omega_{\beta, \alpha}$;  indeed, we have that
\begin{align*}
\ \ \ \ \ \ \ \ \Omega_{\beta, \alpha}^{+} &= \left\{ \omega \in \{ 0, 1\}^{\mathbb{N}} \colon k(0) \preceq \sigma^{n}(\omega) \prec k_{-} \, \ \textup{or} \ \, k_{+} \preceq \sigma^{n}(\omega) \preceq k(1) \, \text{for all} \, n \in \mathbb{N}_{0} \right\}\!,\\
\Omega_{\beta, \alpha}^{-} &= \left\{ \omega \in \{ 0, 1\}^{\mathbb{N}} \colon k(0) \preceq  \sigma^{n}(\omega) \preceq k_{-} \, \ \textup{or} \ \, k_{+} \prec \sigma^{n}(\omega) \preceq k(1) \, \text{for all} \, n \in \mathbb{N}_{0} \right\}\!.
\end{align*}
Here, $\prec$, $\preceq$, $\succ$ and $\succeq$ denote the lexicographic orderings on $\{ 0 ,1\}^{\mathbb{N}}$.  Moreover,
$
\Omega_{\beta, \alpha} = \Omega_{\beta, \alpha}^{+} \cup \Omega_{\beta, \alpha}^{-}$ is closed with respect to the metric $d$ and hence is a subshift.
\end{theorem}
\begin{lemma}[{\cite[Theorem 2]{hubbard1990classification}}]\label{rem:k+k-}
\begin{align*}
\ \ \ \ \ \ \ \ \Omega_{\beta, \alpha}&= \left\{ \omega \in \{ 0, 1\}^{\mathbb{N}} \colon k(0) \preceq \sigma^{n}(\omega) \preceq k_{-} \, \ \textup{or} \ \, k_{+} \preceq \sigma^{n}(\omega) \preceq k(1) \, \ \text{for all} \, n \in \mathbb{N}_{0} \right\}\!,\\
 &= \left\{ \omega \in \{ 0, 1\}^{\mathbb{N}} \colon k(0) \preceq  \sigma^{n}(\omega) \preceq k(1) \ \text{for all} \, n \in \mathbb{N}_{0} \right\}\!.
\end{align*}
\end{lemma}
The connection between intermediate $\beta$-transformations and the $\beta$-expansions of real numbers is given via the $T_{\beta, \alpha}^{\pm}$-expansions of a point. The projection $\pi_{\beta, \alpha} \colon \{ 0, 1 \}^{\mathbb{N}} \to [0, 1]$ is defined as
\begin{align*}
\pi_{\beta, \alpha}(\omega_{1} \omega_{2}\dots) \coloneqq \frac{\alpha}{1 - \beta} + \sum_{k = 1}^{\infty} \frac{\omega_{k}}{\beta^k},
\end{align*}
such that the following diagram commutes.
\begin{align*}
\begin{array}
[c]{ccc}
\Omega_{\beta, \alpha}^{\pm} & \overset{\sigma}{\longrightarrow} & \Omega_{\beta, \alpha}^{\pm}\\
& & \\
\pi_{\beta, \alpha} \downarrow \uparrow \tau^{\pm}_{\beta, \alpha} &  & \pi_{\beta, \alpha} \downarrow \uparrow \tau^{\pm}_{\beta, \alpha} \\ & & \\
 \lbrack 0,1 \rbrack & \underset{T_{\beta, \alpha}^{\pm}} {\longrightarrow} & \lbrack 0,1\rbrack
\end{array}
\end{align*}
\
\par Equivalent conditions, in terms of the kneading invariants $(k_{+},k_{-})$, for the finite type property of $\beta$-shifts are as follows.

\begin{theorem}[{\cite[Theorem 2.3]{li2019denseness}}]\label{lazysft}
For $\beta \in (1, 2)$, we have that
\begin{enumerate}
\item\label{thmA:2} the greedy $\beta$-shift (that is when $\alpha = 0$) is a \textsl{SFT} if and only if $k_{-}$ is periodic;
\item\label{thmA:3} the lazy $\beta$-shift (that is when $\alpha = 2 - \beta$) is a \textsl{SFT} if and only if $k_{+}$ is periodic.
\end{enumerate}
\end{theorem}
\begin{theorem}[{\cite[Theorem 1.3]{li2016intermediate}}]\label{intersft}
Let $(\beta,\alpha) \in \Delta$, the intermediate $\beta$-shift $\Omega_{\beta, \alpha}$ is a \textsl{SFT} if and only if both $k_{+}$ and $k_{-}$  are periodic.
\end{theorem}

Notice that $k_{+}$ and $k_{-}$ are periodic if and only if $k(0)$ and $k(1)$ are both periodic. A necessary and sufficient condition, in terms of the kneading invariants, for the property of an intermediate $\beta$-shift to be a sofic shift can be found in \cite[Proposition 2.14]{kalle2012beta} and stated as, the intermediate $\beta$-shift $\Omega_{\beta, \alpha}$ is sofic if and only if both $k_{+}$ and $k_{-}$  are eventually periodic. The following result shows that kneading invariants are symmetric about the middle point $(\beta, 1-\beta/2)$. Let $(w_{1},w_{2}\cdots)\in\{0,1\}^{\mathbb{N}}$, the \textsl{symmetric map} $s:\{0,1\}^{\mathbb{N}}\circlearrowleft$ is defined by
$$ s(w_{1},w_{2}\cdots)=(w_{1}+1 (\bmod 2), w_{2}+1 (\bmod 2),\cdots).
$$
\begin{lemma}[{\cite[Theorem 2.5]{li2016intermediate}}]\label{symmetric}
For $(\beta,\alpha)\in\Delta$, by the symmetric map $s$, we have $\tau^{+}_{\beta, \alpha}(c_{1})=s(\tau^{-}_{\beta, 2-\beta-\alpha}(c_{2}))$ and $\tau^{-}_{\beta, \alpha}(c_{1})=s(\tau^{+}_{\beta, 2-\beta-\alpha}(c_{2}))$, where $c_{1}=(1-\alpha)/\beta$ and $c_{2}=(\beta+\alpha-1)/\beta$. Moreover, given $\beta\in(1,2)$, there exists a unique point $\alpha=1-\beta/2$, such that $\tau^{+}_{\beta, \alpha}(c)=s(\tau^{-}_{\beta, \alpha}(c))$ and critical point $c=1/2$.
\end{lemma}

\subsection{Combinatorial renormalization}\label{determinant}
\
\par Intermediate $\beta$-transformations are also called linear Lorenz maps. A Lorenz map on $I=[0,1]$ is an interval map $f:I \to I$ such that for some $c\in (0,1)$ we have: (1) $f$ is strictly increasing on $[0,c)$ and on $(c,1]$; (2) $\lim_{x \uparrow c}f(x)=1$, $\lim_{x \downarrow c}f(x)=0$. If, in addition, $f$ satisfies the topological expansive condition: (3) The pre-images set $C(f)=\cup_{n \ge 0}f^{-n}(c)$ of $c$ is dense in $I$, then $f$ is said to be an expansive Lorenz map.
There are lots of studies about Lorenz maps, such as renormalization \cite{cui2015,ding2011,ding2022alpha,glendinning1993prime}, kneading invariants \cite{glendinning1996zeros,hubbard1990classification,glendinning1993prime}, conjugate invariants \cite{cui2015,ding2021complete,glendinning1990topological} and so on. Similar to Section 2.2, given an expansive Lorenz map $f$, we can also obtain its kneading invariants and Lorenz-shift, denoted as $\Omega_{f}$. Now we consider the renormalization of Lorenz map in combinatorial
way. The following definition is essentially from Glendinning and Sparrow\cite{glendinning1993prime}.

\begin{definition}
Let $f$ be an expansive Lorenz map, we say the kneading invariants $(k_{+},k_{-})$ of $\Omega_f$ is \textsl{renormalizable} if there exists finite, non-empty words $(w_+,\ w_-)$, such that
\begin{equation} \label{com-renormal1}
\left \{ \begin{array}{ll}
k_{+}=1k(0) =&w_+ w_-^{p_1} w_+^{p_2} \cdots,\\
k_{-}=0k(1) =&w_- w_+^{m_1} w_-^{m_2} \cdots,
\end{array}
\right.
\end{equation}
where $w_+=1\cdots$, $w_-=0\cdots$, the lengths $|w_+|>1$
and $|w_-|> 1$, and $p_1, m_1>0$. The kneading invariants of the \textsl{renormalization} is $RK=(k_+^1, k_-^1)$, where
\begin{equation} \label{com-renormal2}
\left \{ \begin{array}{ll}
k_{+}^{1}=1k(0)_{1}=&1 0^{p_1} 1^{p_2} \cdots,\\
k_{-}^{1}=0k(1)_{1}= & 0 1^{m_1} 0^{m_2} \cdots.
\end{array}
\right.
\end{equation}
\end{definition}

To describe the renormalization more concisely, we use $*$-product, which is introduced for kneading invariants of Lorenz map in 1990s\cite{glendinning1990topological}. The $*$-product of renormalizaition is defined to be $K=(k_+, k_-):=W*RK$, i.e.,
\begin{equation}\label{product}
(k_{+},k_{-})=(w_+,\ w_-)*(k_{+}^1,k_{-}^1)
\end{equation}
where $(k_{+},k_{-})$ is the pair of sequence obtained by
replacing $1$'s by $w_+$, replacing $0$'s by $w_-$ in $k_{+}^1$ and $k_{-}^1$. Using $*$-product, (\ref{com-renormal1}) and (\ref{com-renormal2}) can be expressed by (\ref{product}). So the kneading invariants are renormalizable if and only if it can be decomposed as the $*$-product; otherwise, we say $(k_{+},k_{-})$ is \textsl{prime}. Note that we do not involve $(w_+, w_-)=(1, 01)$ and $(w_+, w_-)=(10, 0)$ in the definition of renormalization as in \cite{glendinning1993prime}. The cases $(w_+, w_-)=(1, 01)$ and $(w_+, w_-)=(10, 0)$  correspond to trivial renormalizations.
\par A renormalization is called \textsl{periodic renormalization} if the renormalization words $w_{+}$ and $w_{-}$ are related with rational rotations or Farey words; otherwise, it is called \textsl{non-periodic renormlization}. Periodic renormalization is also called primary $n(k)$-cycle in \cite{glendinning1990topological,palmer1979classification}, see the following example for intuitive understanding.

\begin{example}\label{periodic} Fix $RK=((100)^{\infty},(01)^{\infty})$, using $\ast$ product, when $w_{+}=10$, $w_{-}=01$, we have $k_{+}=(100101)^{\infty}$ and $k_{-}=(0110)^{\infty}$, which corresponds to 2(1)-cycle and rational number $1/2$. Similarly, $w_{+}=100$ and $w_{-}=010$ correspond to $1/3$; $w_{+}=101$ and $w_{-}=011$ corresponds to $2/3$.
\end{example}

\subsection{Linearizable kneading invariants}

\
\par We know that each pair of $(\beta,\alpha)\in\Delta$ corresponds to a $\beta$-shift $\Omega_{\beta,\alpha}$, and $\Omega_{\beta,\alpha}$ is totally determined by kneading invariants $(k_{+},k_{-})$. A natural question stemming from above is that, what kind of sequences in $\{ 0, 1 \}^{\mathbb{N}}$ can be the kneading invariants of an intermediate $\beta$-transformation? Actually, this is a question about when an expansive Lorenz map $f$ can be uniformly linearized. We first see the following admissible condition.
\begin{theorem}[{\cite[Theorem 1]{hubbard1990classification}}]\label{kneading space}
If $f$ is a topologically expansive Lorenz map, then the kneading invariants $(k_{+},k_{-})$, satisfies
	\begin{equation} \label{expansive lorenz}
		\sigma (k_{+})\preceq\sigma^n(k_{+})\prec\sigma (k_{-}), \ \ \ \sigma (k_{+})\prec \sigma^n(k_{-})\preceq\sigma (k_{-})  \ \ \ {\rm for\ all} \ \  n \ge 0,
	\end{equation}
 Conversely, given any two sequences $k_{+}$ and $k_{-}$ satisfying $(\ref{expansive lorenz})$, there exists an expansive Lorenz map $f$ with $(k_{+},k_{-})$ as its kneading invariant, and $f$ is unique up to conjugacy.
\end{theorem}
\par We call kneading invariants $(k_{+},k_{-})$ is admissible if condition (2.4) above is satisfied.
When considering the subshift induced by expansive Lorenz map $f$, $\Omega_{f}$ can be prime, periodically renormalized or non-periodically renormalized. However, $f$ can topologically conjugate to a $\beta$-transformation if and only if $f$ is prime or can only be periodically renormalized for finite times \cite{cui2015,glendinning1990topological}, and we call that $f$ can be uniformly linearized.
\begin{lemma}[{\cite[Theorem A]{glendinning1990topological}}{\cite[Theorem A]{cui2015}}]\label{finite-periodic}
When $\sqrt{2}<\beta<2$, for all $\alpha\in(0,2-\beta)$, $\Omega_{\beta,\alpha}$ is prime. When  $1<\beta\leq\sqrt{2}$ and $\alpha\in(0,2-\beta)$, $\Omega_{\beta,\alpha}$ is either prime or can only be periodically renormalized for finitely many times.
\end{lemma}
Combining with Theorem \ref{kneading space} and Lemma \ref{finite-periodic}, we are able to give the definition of linearizable kneading invariants for expansive Lorenz maps $f$.
Since here we consider the $\beta$-shifts with $\beta>1$, then extra condition is needed to make sure the topological entropy $h_{top}(\sigma, \Omega_{\beta,\alpha})=\log\beta>0$. Similar to Remark \ref{rem:k+k-}, denote $$\Omega(k_{+}, k_{-}):=\{\omega\in\{0,1\}^{\mathbb{N}}: \sigma(k_{+})\preceq\sigma^{n}(\omega) \preceq \sigma(k_{-}) \ \text{for all} \, n \in \mathbb{N}_{0}\}.$$
\begin{definition}\label{defn:admissible}
A pair of infinite sequences $(k_{+},k_{-})$ is said to be \textsl{linearizable} to a $\beta$-shift with $\beta>1$ if the following conditions are satisfied:
\begin{enumerate}
\item\label{enumi2:defn_admissible} 
$(k_{+},k_{-})$ is admissible,
\item\label{enumi4:defn_admissible} $(k_{+},k_{-})$ is prime or can only be periodically renormalized for finite times.

\item\label{enumi3:defn_admissible} $\displaystyle{\lim_{n \to \infty} n^{-1} \ln \left( \# \Omega(k_{+},k_{-})\lvert_{n}  \right) > 0}$,
\end{enumerate}
If in addition $k_{+}$ and $k_{-}$ are periodic, then we call $(k_{+},k_{-})$ is \textsl{periodically linerizable}.
\end{definition}

With the help of Definition \ref{defn:admissible}, we are able to answer the question at the initial of this section. The limit in \Cref{defn:admissible} \eqref{enumi3:defn_admissible} exists for the sequence $\{ \ln(\# \Omega(k_{+},k_{-})\vert_{n} \}_{n \in \mathbb{N}}$ is sub-additive. In fact, the limit is the topological entropy of the kneading space induced by $(k_{+},k_{-})$. 

\begin{remark}
Two infinite sequences $k_{+},k_{-}\in \{ 0, 1\}^{\mathbb{N}}$ are kneading invariants for an intermediate $\beta$-shift if and only if $(k_{+},k_{-})$ is a linearizable pair.
\end{remark}
\begin{corollary}\label{periodic admissible}
Two infinite sequences $k_{+},k_{-} \in \{ 0, 1\}^{\mathbb{N}}$ are kneading invariants for an intermediate $\beta$-shift of finite type if and only if $(k_{+},k_{-})$ is a periodically linearizable pair.
\end{corollary}

\subsection{Kneading determinant}\label{determinant}
\
\par The idea for kneading determinant goes back to \cite{milnor1988iterated}; see also \cite{glendinning1996zeros}. Let $({k_ + },{k_ - })$ be the kneading invariants of an expansive Lorenz map $f$, where $k_{+}=(v_{1}v_{2}\cdots)$ and $k_{-}=(w_{1}w_{2}\cdots)$. Then the kneading determinant is a formal power series defined as $K(t) = {K_ + }(t) - {K_ - }(t)$ , where
$${K_ + }(t) = \sum\limits_{i = 1}^\infty  {v_{i}{t^{i-1}}},\ \ \ \ {K_ - }(t) = \sum\limits_{i = 1}^\infty  {w_{i}{t^{i-1}}}.$$
 The following lemma offers a straight method to calculate topological entropy of $\Omega_{f}$ if the kneading invariants are given.
\begin{lemma}[{\cite[Lemma 3]{barnsley2014} \cite[Theorem 3]{glendinning1996zeros}}]\label{zeros}
Let $({k_ + },{k_ - })$ be a pair of infinite sequences satisfying conditions $(1)$, $(3)$ in Definition \ref{defn:admissible}. Let $K(t)$ be the corresponding kneading determinant, and $t_{0}$ be the smallest positive root of $K(t)$ in $(0,1)$, then  we have $h_{top}(\sigma,\Omega(k_{+},k_{-}))=-\log t_{0}$.
\end{lemma}

\begin{remark}
If $(k_{+},k_{-})$ is a pair of linearizable kneading invariants, i.e, $(k_{+},k_{-})$ corresponds to intermediate $\beta$-shift, then $1/\beta$ equals the smallest positive root of $K(t)$ in $(0,1)$.

\end{remark}

Two linearizable kneading invariants with the same numerator of kneading determinants will lead to the same smallest positive root, therefore they have the identical $\beta$. However, a fixed $\beta$ may come from different kneading determinants. For instance, $(k_{+},k_{-})=((1000)^{\infty},(01)^{\infty})$ and $(k^{\prime}_{+},k^{\prime}_{-})=((1001)^{\infty},(011)^{\infty})$ have the same $\beta\approx1.46557$, but different kneading determinants.

Given a pair of linearizable kneading invariants $(k_{+},k_{-})$, we know that $\beta=1/t_{0}$, where $t_{0}$ is the smallest positive root of $K(t)$, but how about $\alpha$? There are several identical formulas of $\alpha$ mentioned in \cite{ding2021complete}, for convenience, here we use $\alpha  = (1/t_{0}  - 1)(K_{+}(t_{0})-1)$. Denote $\mathbb{Q}(\beta)$ as the smallest sub-field of the reals containing $\beta$, by the formula of $\alpha$, we now that if $(\beta,\alpha)\in\mathcal{F}(\beta)$, then $\alpha\in\mathbb{Q}(\beta)$. Hence $\mathcal{F}(\beta)$ is also a countable set. The following lemma shows that, when \Cref{defn:admissible} \eqref{enumi4:defn_admissible} is violated, different pairs of infinite sequences may be figured out the same parameter $(\beta,\alpha)$, but only one of them can be uniformly linearized to a $\beta$-shift.

\begin{lemma}[{\cite[Lemma 8]{barnsley2014}}{\cite[Proposition 2]{ding2021complete}}]\label{xunibeta}
Let $k_{+} $ and $k_ {- }$ be infinite sequences with $k_{+}|_{2}=10$, $k_{-}|_{2}=01$, and satisfying condition Definition \ref{defn:admissible} $(1)$, $(3)$.  Let $(\beta,\alpha)$ be the corresponding parameter pair, and $c = (1 - \alpha)/\beta$.
\par (1) If $\tau^{+}_{\beta, \alpha}(c)=k_{+}$ and $\tau^{-}_{\beta, \alpha}(c)=k_{-}$, then $(k_{+},k_{-})$ satisfies Definition \ref{defn:admissible} and $(k_{+},k_{-})$ is the linearizable kneading invariants of $\Omega_{\beta,\alpha}$.
\par (2) If $\tau^{+}_{\beta, \alpha}(c)\neq k_{+}$ or $\tau^{-}_{\beta, \alpha}(c)\neq k_{-}$, then $(k_{+},k_{-})$ can be non-periodic renormalized. Suppose the $m$-th $(1\leq m<\infty)$ renormalization is the nearest non-periodic renormalization with words $(w_{+},w_{-})$, and the $i$-th periodic renomalization words are $(w_{+,i},w_{-,i})$ $(0\leq i\leq m-1)$, we have $(\tau^{+}_{\beta, \alpha}(c),\tau^{-}_{\beta, \alpha}(c))=(w_{+,1},w_{-,1})\ast\cdots\ast(w_{+,m-1},w_{-,m-1})\ast(w^{\infty}_{+},w^{\infty}_{-})$.
\end{lemma}
In other words, if condition $(2)$ of Definition \ref{defn:admissible} breaks, then $(k_{+},k_{-})$ will not correspond to an intermediate $\beta$-transformation, but it still can be calculated and obtain a pair of virtual parameter $(\beta,\alpha)$, Lemma \ref{xunibeta} shows how to obtain the linearizable kneading invariants with respect to $(\beta,\alpha)$. See Example \ref{realkneading} for an intuitive explanation.

\section{Proof of main results}
\par We show that, for each $(\beta,\alpha)\in\Delta$ with $\Omega_{\beta,\alpha}$ being a SFT, the  intermediate $\beta$-transformation $T_{\beta,\alpha}$ has matching. Moreover, for any $(\beta,\alpha)\in\mathcal{M}$, there exist $(\beta,\alpha^{\prime})\in \Delta(\beta)$ and $\alpha^{\prime}$ is $\epsilon$-close to $\alpha$ under the Euclidean metric, such that $\Omega_{\beta,\alpha^{\prime}}$ is a SFT. We also consider the endpoints of matching interval $I(\beta,\alpha)$ when $(\beta,\alpha)\in\mathcal{M}$, and give a classification of the endpoints. Finally, we prove Corollary \ref{cor1.3}. \\

\noindent \textbf{Proof of Theorem \ref{th1} (1)}
\par Let $(\beta,\alpha)\in\mathcal{F}$, we need to prove $(\beta,\alpha)\in\mathcal{M}$.
By Corollary \ref{periodic admissible}, the kneading invariants for an intermediate $\beta$-shift of finite type are both periodic. Hence we suppose
$$
\left \{ \begin{array}{ll}
k_{+}=(10a_{3}a_{4}\cdots a_{n})^{\infty},\\
k_{-}=(01b_{3}b_{4}\cdots b_{n})^{\infty}.
\end{array}
\right.
$$
Notice that if the periods are different, take $n$ as their least common multiple. According to the symbols of $k_{+}$, we have $(T^{+}_{\beta,\alpha})^{1}(c)=0$, $(T^{+}_{\beta,\alpha})^{2}(c)=\alpha$, $(T^{+}_{\beta,\alpha})^{3}(c)=\beta\alpha+\alpha-a_{3}$,
$(T^{+}_{\beta,\alpha})^{4}(c)=\beta^2\alpha+\beta\alpha-\beta a_{3}+\alpha-a_{4}$, $(T^{+}_{\beta,\alpha})^{5}(c)=\beta^3\alpha+\beta^2\alpha-\beta^2 a_{3}+\beta\alpha-\beta a_{4}+\alpha-a_{5}$.
By induction,
$$ (T^{+}_{\beta,\alpha})^{n}(c)=\beta^{n-2}\alpha+\beta^{n-3}\alpha-\beta^{n-3} a_{3}+\cdots +\beta\alpha-\beta a_{n-1}+\alpha-a_{n}.
$$
Similarly, with the coding of $k_{-}$, we have
$$ (T^{-}_{\beta,\alpha})^{n}(c)=\beta^{n-1}+\beta^{n-2}\alpha-\beta^{n-2}+\beta^{n-3}\alpha-\beta^{n-3}b_{3}+\cdots +\beta\alpha+\beta b_{n-1}+\alpha-b_{n}.
$$
\par By Lemma \ref{zeros}, the kneading determinant of $(k_{+},k_{-})$  is
$$K(t)=(1-t+(a_{3}-b_{3})t^2+\cdots +(a_{n}-b_{n})t^{n-1})/(1-t^{n}),$$
and $1/\beta$ is the smallest positive root in $(0,1)$. Replacing $t$ with $1/\beta$, we denote
$$ D(\beta):=\beta^{n-1}-\beta^{n-2}+(a_{3}-b_{3})\beta^{n-3}+\cdots+(a_{n}-b_{n})=0.
$$
It can be seen that $(T^{-}_{\beta,\alpha})^{n}(c)-(T^{+}_{\beta,\alpha})^{n}(c)=D(\beta)$=0. However, $k_{+}$ and $k_{-}$ may have the same tail, which means there exists a minimal integer $m\leq n$ such that $a_{i}=b_{i}$ for $m<i\leq n$. In this way, we have
\begin{align*}
 D(\beta)&=  \beta^{n-1}-\beta^{n-2}+(a_{3}-b_{3})\beta^{n-3}+\cdots+(a_{m}-b_{m})\beta^{n-m}\!\\
 &=\beta^{m-1}-\beta^{m-2}+(a_{3}-b_{3})\beta^{m-3}+\cdots+(a_{m}-b_{m})\!\\
 &=(T^{-}_{\beta,\alpha})^{m}(c)-(T^{+}_{\beta,\alpha})^{m}(c)=0.\!
\end{align*}
This indicates  $m$ is the minimal integer such that $(T^{+}_{\beta,\alpha})^{m}(c)=(T^{-}_{\beta,\alpha})^{m}(c)$. Since $T^{+}_{\beta,\alpha}(c)=0$ and $T^{-}_{\beta,\alpha}(c)=1$, we have $(T^{+}_{\beta,\alpha})^{m-1}(0)=(T^{-}_{\beta,\alpha})^{m-1}(1)$, which means matching occurs at the time $m-1$, and $(\beta,\alpha)\in\mathcal{M}$. Then $\mathcal{F}\subset\mathcal{M}$. However, there exist examples that $(\beta,\alpha)\in\mathcal{M}$, but $(\beta,\alpha)\notin\mathcal{F}$, see Example \ref{example1}. Hence we have $\mathcal{F}\subsetneq\mathcal{M}$.
$\hfill\square$

\begin{proposition}\label{SFTIFF}  Let $(\beta,\alpha)\in \mathcal{M}$. Then $(\beta,\alpha)\in\mathcal{F}$ if and only if $0$ and $1$ iterate back to the critical point again.

\begin{proof}
 Suppose $T_{\beta,\alpha}$ has matching at time $m-1$, then $(T^{+}_{\beta,\alpha})^{m}(c)=(T^{-}_{\beta,\alpha})^{m}(c)$. Let $k_{+}=(10\cdots a_{m}\cdots)$ and $k_{-}=(01\cdots b_{m}\cdots)$,
there are two cases about the tail of $k_{+}$ and $k_{-}$.
\par Case one, the iterations of $0$ and $1$ never get back to the critical point again, which means $(T^{+}_{\beta,\alpha})^{m+i}(c)=(T^{-}_{\beta,\alpha})^{m+i}(c)\neq c$ for any $i\geq0$, and
$$\left \{ \begin{array}{ll}
k_{+}=10\cdots a_{m}w_{m+1}w_{m+2}\cdots,\\
k_{-}=01\cdots b_{m}w_{m+1}w_{m+2}\cdots.
\end{array}
\right.
$$
At this case, we have $a_{i}=b_{i}$ for all $i>m$, and neither $k_{+}$ nor $k_{-}$ is periodic. Notice that $k_{+}$ and $k_{-}$ can be both eventually periodic.
\par Case two, the iterations of $0$ and $1$ get back to the critical point again, that is, there exists smallest integers $i$ and $j$ such that $(T^{+}_{\beta,\alpha})^{m+i}(c)=(T^{-}_{\beta,\alpha})^{m+j}(c)= c$. At this case,
$$\left \{ \begin{array}{ll}
k_{+}=(10\cdots a_{m}w_{m+1}\cdots w_{m+i})^{\infty},\\
k_{-}=(01\cdots b_{m}w_{m+1}\cdots w_{m+j})^{\infty},
\end{array}
\right.
$$
and $(\beta,\alpha)\in\mathcal{F}$. It can be seen that both cases will lead to the same $\beta$, which indicates that the essential information of $\beta$ depends on the first $m$ symbols of $(k_{+},k_{-})$ if matching occurs at time $m-1$.
\end{proof}
\end{proposition}
\begin{lemma}\label{denghaobudui}
 Let $k_{+}=(v_{1}\cdots v_{k}\cdots)$, $k_{-}=(w_{1}\cdots w_{k}\cdots)$, and $v_{i}=w_{i}$ for all $i>k$. Then $(k_{+},k_{-})$ is not linearizable if the following equality holds for some $r<k$,
\begin{equation}\label{equation}
(v_{1}\cdots v_{k})=(w_{r+1}\cdots w_{k}w_{k+1}\cdots w_{k+r})=(w_{r+1}\cdots w_{k}w_{1}\cdots w_{r}).
\end{equation}

\begin{proof}
\par Since $v_{i}=w_{i}$ for all $i>k$, here denote $v=10v_{3}\cdots v_{k}\mu$, $w=01w_{3}\cdots w_{k}\mu$, where $\mu$ is the infinite sequence $(w_{k+1}w_{k+2}\cdots)$. 
By equality (\ref{equation}), we have $(w_{k+1}\cdots w_{k+r})=(w_{1}\cdots w_{r})$, $w_{r+1}=1$ and $w_{k+1}=0$. Then we consider $w_{k+r+1}$ and divide it into two cases.
\par Case one, $w_{k+r+1}=0$. At this case, we compare $(w_{k+r+1}\cdots w_{k+2r})$ with $(w_{1}\cdots w_{r})$. If $(w_{k+r+1}\cdots w_{k+2r})\succ(w_{1}\cdots w_{r})$, then $\sigma^{k+r}(k_{-})\succ k_{-}$, $(k_{+},k_{-})$ violates Definition \ref{defn:admissible} (1); if $(w_{k+r+1}\cdots w_{k+2r})\prec(w_{1}\cdots w_{r})$, then $(w_{r+1}\cdots w_{k+r}w_{k+r+1}\cdots w_{k+2r})\prec(v_{1}\cdots v_{k}v_{k+1}\cdots v_{k+r})$, which implies $\sigma^{r}(k_{-})\prec k_{+}$ and Definition \ref{defn:admissible} (1) is violated.
\par Case two, $w_{k+r+1}=1$. At this case, we compare $(w_{k+r+1}\cdots w_{2k})$ with $(v_{1}\cdots v_{k-r})$. Similarly, when $w_{k+r+1}=1$, both $(w_{k+r+1}\cdots w_{2k})\succ(v_{1}\cdots v_{k-r})$ (leads to $k_{-}\prec\sigma^{k}(k_{-})$) and $(w_{k+r+1}\cdots w_{2k})\prec(v_{1}\cdots v_{k-r})$ (leads to $k_{+}\succ\sigma^{k+r}(k_{-})$) will also indicate that Definition \ref{defn:admissible} (1) is violated.
\par Continue the above process, unless $\mu$ purely consists of $(w_{1}\cdots w_{r})$ and $(v_{1}\cdots v_{k-r})$, we can always find an integer $l\in\mathbb{N}$ with $(w_{l+r+1}\cdots w_{l+2r})\succ ({\rm or} \prec) (w_{1}\cdots w_{r})$ or $(w_{l+r+1}\cdots w_{l+k})$ $\succ ({\rm or} \prec) (v_{1}\cdots v_{k-r}))$ such that $(k_{+},k_{-})$ is not admissible. And at the case $\mu$ purely consists of $(w_{1}\cdots w_{r})$ and $(v_{1}\cdots v_{k-r})$, $(k_{+},k_{-})$ can be renormalized via $w_{+}=(v_{1}\cdots v_{k-r})$ and $w_{-}=(w_{1}\cdots w_{r})$. If $(w_{+},w_{-})$ is a pair of non-periodic renormalization words, then $(k_{+},k_{-})$ violates Definition \ref{defn:admissible} (2); if $(w_{+},w_{-})$ is a pair of periodic renormalization words, then $(k_{+},k_{-})$ violates condition (3) in Definition \ref{defn:admissible}.

\end{proof}

\end{lemma}

\begin{proposition} \label{bijin}
For any $\epsilon>0$ and $(\beta,\alpha)\in\mathcal{M}$, there exists $\alpha^{\prime}$ satisfying $|\alpha-\alpha^{\prime}|<\epsilon$, such that $(\beta,\alpha^{\prime})\in\mathcal{F}(\beta)$.
\begin{proof}
\par Suoppose $(\beta,\alpha)\in\mathcal{M}$ and matching occurs at time $m-1$ ($m\geq2$ be an integer), and $(v=\tau^{+}_{\beta,\alpha}(c),w=\tau^{-}_{\beta,\alpha}(c))$ be the corresponding kneading invariants, where $c=\frac{1-\alpha}{\beta}$. By the admissibility of $(v, w)$, we have
\begin{equation}\label{yunxutj}
		\sigma(v)\preceq\sigma^{s}(v)\prec\sigma(w), \ \ \ \ \sigma(v)\prec\sigma^{s}(w)\preceq\sigma(w), \ \ \ \forall s \ge 0.
\end{equation}
If $(\beta,\alpha)\in\mathcal{F}$, then $(\beta,\alpha)$ can be approached by itself, the proof is completed. If $(\beta,\alpha)\in\mathcal{M}\setminus\mathcal{F}$, by Corollary \ref{periodic admissible}, both $v$ and $w$ are not periodic. Let $v=(v_{1}v_{2}\cdots v_{m}\cdots)$ and $w=(w_{1}w_{2}\cdots w_{m}\cdots)$, by Proposition \ref{SFTIFF}, $v_{i}=w_{i}$ for all $i>m$. Our aim is to construct a pair of periodically linearizable kneading invariants $(k_{+}=\tau^{+}_{\beta,\alpha^{\prime}}(c^{\prime}),k_{-}=\tau^{-}_{\beta,\alpha^{\prime}}(c^{\prime}))$, arbitrarily close to the pair $(v,w)$, where $c^{\prime}=\frac{1-\alpha^{\prime}}{\beta}$ and $(\beta,\alpha^{\prime})\in\mathcal{F}(\beta)$. To complete the proof we also show that the Euclidean distance between $\alpha^{\prime}$ and $\alpha$ is arbitrarily small.
\par \textbf{Step 1. The Construction of $(k_{+},k_{-})$.} Without loss of generality, here we start the construction with $w$, and it is similar to start with $v$. Since $\alpha\neq 2-\beta$, then $w\neq 01^\infty$, and we may choose an integer $n> m$ with $w_{n}=0$. Let $ j\geq 1$ be the minimal integer such that $(w_{j+1}\cdots w_{n})=(w_{1}\cdots w_{n-j})$; note this equality holds for $j=n-1$. Since $w$ is not periodic, we have $\sigma^{j}(w)\prec w$. Therefore, there exists a minimal integer $k\geq n$ with $w_{k+1}=0$ and $w_{k-j+1}=1$, in which case,
\begin{equation}\label{equal}
(w_{j+1}\cdots w_{n-1}w_{n}\cdots w_{k})=(w_{1}\cdots w_{n-j-1}w_{n-j}\cdots w_{k-j}).
\end{equation}
\par Notice that, given an $n>m$ with $w_{n}=0$, there exist corresponding $j\leq n-1$ and $k\geq n$, and such $n$ can be infinitely many. As $n$ varies, $j$ and $k$ also change, and $k$ can also be infinitely many, denote $\mathcal{K}$ be the collection of all such $k$. Let $k_{-}=(w_{1}\cdots w_{m}w_{m+1}\cdots w_{k})^\infty$, and correspondingly
$k_{+}=(v_{1}\cdots v_{m}w_{m+1}\cdots w_{k})^\infty$. Since $w_{k+1}=0$ and $w$ is not periodic, then $\sigma^{k}(w)\prec w$, i.e.,
$$w_{k+1}\cdots w_{2k}w_{2k+1}\cdots w_{3k}\cdots \prec w_{1}\cdots w_{k}w_{k+1}\cdots w_{2k}\cdots=w.
$$
Combing this with the fact $k_{-}=(w_{1}\cdots w_{k})^{\infty}$, we conclude that $w\prec k_{-}$. On the other hand, for the reason that $w_{k+1}=0$ and $v_{1}=1$, we have $v=v_{1}\cdots v_{m}w_{m+1}\cdots w_{k}w_{k+1}\cdots\prec(v_{1}\cdots v_{m}w_{m+1}\cdots w_{k})^{\infty}=k_{+}$. Our aim is to prove $k_{+}$ and $k_{-}$ satisfying
\begin{equation}\label{aim}
\sigma(k_{+})\preceq\sigma^{r}(k_{+})\prec\sigma(k_{-}), \ \ \sigma(k_{+})\prec\sigma^{r}(k_{-})\preceq\sigma(k_{-}), \ \ \forall r \in\{0,\cdots,k-1\}.
\end{equation}

\par \textbf{Step 2. The linearization of $(k_{+},k_{-})$.} First, we claim that $\sigma^{r}(k_{-})\preceq\sigma(k_{-})$ for all $r \in\{0,\cdots,k-1\}$. When $w_{r+1}=0$, $\sigma^{r}(k_{-})\prec\sigma(k_{-})$ is obvious. It remains to show that $\sigma^{r}(k_{-})\preceq k_{-}$ for all $r\in \{0,\cdots,k-1\}$ with $w_{r+1}=0$, as this implies that $\sigma^{r}(k_{-})\preceq\sigma (k_{-})$ with $w_{r+1}=1$. Let $w_{r+1}=0$. When $0<r<j$, $\sigma^{r}(w)\preceq w$ and the minimality of $j$ indicates $(w_{r+1}\cdots w_{n})\prec(w_{1}\cdots w_{n-r})$, hence $\sigma^{r}(k_{-})\prec k_{-}$; when $j\leq r\leq k-1$, we have
\begin{align*}
(w_{r+1}\cdots w_{k} w_{1})&=(w_{r+1}\cdots w_{k} w_{k+1})\!\\
 &\prec(w_{r+1}\cdots w_{k} 1)\!\\
 &=(w_{r-j+1}\cdots w_{k-j} w_{k-j+1})\!\\
  &\preceq(w_{1}\cdots w_{k-r} w_{k-r+1}),\!
\end{align*}
which yields $\sigma^{r}(k_{-})\prec k_{-}$. When $r=1$, we have $\sigma^{r}(k_{-})=\sigma(k_{-})$. Therefore, $\sigma^{r}(k_{-})\preceq\sigma(k_{-})$ for all $r \in\{0,\cdots,k-1\}$.
\par Second, we prove $\sigma(k_{+})\preceq\sigma^{r}(k_{+})$ for any $r\in\{0,\cdots, k-1\}$, where $k_{+}=(v_{1}\cdots v_{m}$ $v_{m+1}\cdots v_{k})^{\infty}$ and $v_{i}=w_{i}$ for $i\in\{m+1,\cdots k\}$. If $v_{r+1}=1$, $\sigma(k_{+})\prec\sigma^{r}(k_{+})$ is clear. It remains to show that $k_{+}\preceq\sigma^{r}(k_{+})$ with $v_{r+1}=1$. When $v_{r+1}=1$, $\sigma^{r}(v)\succeq v$ holds. For any $1< r\leq k-1$, and $v_{k+1}=w_{k+1}=0$, we have
\begin{align*}
(v_{r+1}\cdots v_{k} v_{1})&\succ(v_{r+1}\cdots v_{k} v_{k+1})\succeq(v_{1}\cdots v_{k-r+1}),\!
\end{align*}
which means $\sigma(k_{+})\prec\sigma^{r}(k_{+})$. When $r=1$, $\sigma^{r}(k_{-})=\sigma(k_{-})$. Hence $\sigma(k_{+})\preceq\sigma^{r}(k_{+})$.
\par Third, we prove $\sigma(k_{+})\prec\sigma^{r}(k_{-})$ for all $r\in\{0,\cdots ,k-1\}$. If $w_{r+1}=1$, $\sigma(k_{+})\prec\sigma^{r}(k_{-})$ holds clearly. It remains to show that $k_{+}\prec\sigma^{r}(k_{-})$ for any $r\in \{0,\cdots,k-1\}$ with $w_{r+1}=1$.
Since $w\prec k_{-}$ and their first $k$ symbols are identical, for any $0< r\leq k-1$ with $w_{r+1}=1$, we have $ v\prec \sigma^{r}(w)\prec\sigma^{r}(k_{-}) $, that is
\begin{equation}
(v_{1}\cdots v_{k})\preceq(w_{r+1}\cdots w_{k}w_{k+1}\cdots w_{k+r})\preceq(w_{r+1}\cdots w_{k}w_{1}\cdots w_{r})
\end{equation}
By Lemma \ref{denghaobudui}, $(v,w)$ is linearizable implies that at least one of the $``\preceq"$ above turns into $``\prec"$. Hence we have $(v_{1}\cdots v_{k})\prec(w_{r+1}\cdots w_{k}w_{1}\cdots w_{r})$ and $\sigma(k_{+})\prec\sigma^{r}(k_{-})$ for all $r\in\{0,\cdots k-1\}$.
\par Finally, we intend to prove $\sigma^{r}(k_{+})\prec\sigma(k_{-})$ for all $r\in\{0,\cdots,k-1\}$. If $v_{r+1}=0$, then $\sigma^{r}(k_{+})\prec\sigma(k_{-})$. It remains to show that $\sigma^{r}(k_{+})\prec k_{-}$ for all $r$ with $v_{r+1}=0$.
Observe that some $k\in\mathcal{K}$  may not satisfy $\sigma^{r}(k_{+})\prec k_{-}$ with $v_{r+1}=0$, see Example \ref{manyk}. Denote $\mathcal{K}^{\prime}:=\{k\in\mathcal{K}: \sigma^{r}(k_{+})\prec k_{-} \ {\rm for \ all} \ r \ {\rm with} \ v_{r+1}=0\}$.
\par {\bf Claim :}  $\mathcal{K}^{\prime}\neq\emptyset$.
\par We prove this claim by contradiction. If for any $k\in\mathcal{K}$, there exists a $r_{k}\in\{0,\cdots,k-1\}$ with $v_{r_{k}+1}=0$ such that $\sigma^{r_{k}}(k_{+})\succeq k_{-}$, then
\begin{equation}\label{finally}
(v_{r_{k}+1}\cdots v_{k}v_{1}\cdots v_{r_{k}})\succeq(w_{1}\cdots w_{k-r_{k}}w_{k-r_{k}+1}\cdots w_{k}).
\end{equation}
Since $v_{r_{k}+1}=0$, we have $\sigma^{r_{k}}(v)\prec w$ and $(v_{r_{k}+1}\cdots v_{k})\preceq(w_{1}\cdots w_{k-r_{k}})$. Combining this with (\ref{finally}), we obtain that $(v_{r_{k}+1}\cdots v_{k})=(w_{1}\cdots w_{k-r_{k}})$, and $(v_{1}\cdots v_{r_{k}})\succeq(w_{k-r_{k}+1}\cdots w_{k})$. For the case $(v_{1}\cdots v_{r_{k}})\succ(w_{k-r_{k}+1}\cdots w_{k})$ with $w_{k-r_{k}+1}=1$, we have $v\succ\sigma^{k-r_{k}}(w)$,
and the condition (\ref{yunxutj}) is violated. Next we consider another two cases.

\par Case A: $(v_{1}\cdots v_{r_{k}})\succ(w_{k-r_{k}+1}\cdots w_{k})$ with $w_{k-r_{k}+1}=0$. Here we divide the value of $r_{k}$ into 3 cases. If $r_{k}=j_{k}$, by equality (\ref{equal}), $w_{k-r_{k}+1}=w_{k-j_{k}+1}=1$, which contradicts with the assumption of Case A. If $r_{k}<j_{k}$, since $j_{k}$ is the minimal integer such that equality (\ref{equal}) holds, hence we have $(w_{r_{k}+1}\cdots w_{k})\prec( w_{1}\cdots w_{k-r_{k}})$. It can be verified that 
when $k$ is large enough, we have $r_{k}\geq m$, and then $(v_{r_{k}+1}\cdots v_{k})=(w_{r_{k}+1}\cdots w_{k})\prec( w_{1}\cdots w_{k-r_{k}})$, which conflicts with $(v_{r_{k}+1}\cdots v_{k})=( w_{1}\cdots w_{k-r_{k}})$. It remains to prove the case $j_{k}<r_{k}<k-1$. Since $j_{k}<r_{k}$, by equality (\ref{equal}), we have $(v_{r_{k}+1}\cdots v_{k})=(w_{r_{k}+1}\cdots w_{k})=( w_{1}\cdots w_{k-r_{k}})=(w_{r_{k}-j_{k}+1}\cdots w_{k-j_{k}})$. See Figure 3 for a clearer illustration. However, $w_{k-r_{k}+1}=0$ but $w_{k-j_{k}+1}=1$, this indicates $\sigma^{r_{k}-j_{k}}(w)\succ w$, which violates the condition (\ref{yunxutj}).
\par  Case B: If for some $k\in\mathcal{K}$, there exists $r_{k}$ such that $(v_{1}\cdots v_{r_{k}})=(w_{k-r_{k}+1}\cdots w_{k})$.
Let $s=k-r_{k}$, and $\xi$ be the first $s$ symbols of $w$, that is, $\xi=(w_{1}\cdots w_{k-r_{k}})$.
Let $t=k-a[k/s]$, and $\eta$ be the first $t$ symbols of $v$, where $[k/s]$ means the integer part of $k/s$. With the help of Figure 4, we obtain that $k_{+}=(\eta \xi\cdots \xi)^{\infty}$, $k_{-}=(\xi\eta \xi\cdots \xi)^{\infty}$, and $(k_{+},k_{-})$ can be renormalized via $\eta$ and $\xi$. Similar to the proof of Lemma \ref{denghaobudui}, 
we have $v=\eta \xi^{\infty}$ and $w=\xi\eta \xi^{\infty}$, which violates the linearizability of $(v,w)$ in Definition \ref{defn:admissible}.

\par Next we show that the cardinality of $\mathcal{K}^{\prime}$ is infinite. In fact, if $k\in \mathcal{K}^{\prime}$, then for any $k^{\prime}\in \mathcal{K}$ with $k^{\prime}>k$, $k^{\prime}\in \mathcal{K}^{\prime}$. We prove this by contradiction. Given any $k\in \mathcal{K}^{\prime}$, if $k^{\prime}>k$ and $k^{\prime}\in \mathcal{K}$, but $k^{\prime}\notin \mathcal{K}^{\prime}$, then there exists $r<k^{\prime}$ with $v_{r+1}=0$ such that $\sigma^{r}(k^{\prime}_{+})\succeq k^{\prime}_{-}$, where $(k^{\prime}_{+}, k^{\prime}_{-})$ means the constructed kneading invariants corresponding to $k^{\prime}$.
Then we have
$
(v_{r+1}\cdots v_{k}\cdots v_{k^{\prime}}v_{1}\cdots v_{r})\succeq(w_{1}\cdots w_{k-r}\cdots w_{k}\cdots w_{k^{\prime}}),
$
which implies $(v_{r+1}\cdots v_{k})\succeq(w_{1}\cdots w_{k-r})$. If $(v_{r+1}\cdots v_{k})\succ(w_{1}\cdots w_{k-r})$, it conflicts with $k\in\mathcal{K}^{\prime}$. Hence $(v_{r+1}\cdots v_{k})=(w_{1}\cdots w_{k-r})$ , and $v_{k+1}=w_{k-r+1}=0$. However, $k\in\mathcal{K}^{\prime}$ implies $\sigma^{r}(k_{+})\prec k_{-}$ and
$
(v_{r+1}\cdots v_{k}v_{1}\cdots v_{r})\prec(w_{1}\cdots w_{k-r}w_{k-r+1}\cdots w_{k})$. Hence we obtain $w_{k-r+1}=v_{1}=1$, which contradicts with $v_{k+1}=w_{k-r+1}=0$. This means such $r$ does not exist, $k^{\prime}$ is also in $\mathcal{K}^{\prime}$.

\par In conclusion, we can find infinitely many $k$ such that $(k_{+},k_{-})$ are periodic kneading invariants satisfying Definition \ref{defn:admissible} (1). Since $n$ can be chosen arbitrarily large, then $(k_{+},k_{-})$ can be arbitrarily close to $(v,w)$ with respect to the usual metric. By construction we know that $(k_{+},k_{-})$ and $(v,w)$ have the same numerator of kneading determinants, hence they have the same $\beta$ and $(k_{+},k_{-})$ satisfies Definition \ref{defn:admissible} (3).
By the proof of Lemma \ref{denghaobudui}, we have $(v,w)$ is prime if and only if $(k_{+},k_{-})$ is prime.
And both $(v,w)$ and $(k_{+},k_{-})$ have the same renormalization words, hence $(k_{+},k_{-})$ can only be periodically renormalized. Moreover, the periods of $k_{+}$ and $k_{-}$ are finite, hence $(k_{+},k_{-})$ can be periodically renormalized for finite times. As a result, the $(k_{+},k_{-})$ we constructed satisfies the  Definition \ref{defn:admissible}.
\par \par \textbf{Step 3. Arbitrarily close to $\alpha$.}By the formula $\alpha  = (\beta  - 1)(K_{+}(1/\beta)-1)$, we have
\begin{align*}
|\alpha-\alpha^{\prime}|&= \Big|(\beta-1)\Big(\sum^{\infty}_{i=1}\frac{v_{i}}{\beta^{i-1}}-1-
\big(\sum^{k}_{i=1}\frac{v_{i}}{\beta^{i-1}}\big)/(1-\beta^{-k})+1\Big)\Big|\!\\
 &=(\beta-1)\Big|\sum^{\infty}_{i=1}\frac{v_{i}}{\beta^{i-1}}-
\big(\sum^{k}_{i=1}\frac{v_{i}}{\beta^{i-1}}\big)/(1-\beta^{-k})\Big|\!\\
  &<(\beta-1)\frac{\beta^{-k}}{1-1/\beta}=\beta^{-(k-1)}.\!
\end{align*}
This implies that, for any $\varepsilon>0$, there exists $k\geq n$ such that both $k_{+}$ and $k_{-}$ are periodic, $|\alpha-\alpha^{\prime}|<\beta^{-k}<\varepsilon$, and $(\beta,\alpha^{\prime})\in\mathcal{F}(\beta)$.
\end{proof}

\end{proposition}

\par As a result of Proposition \ref{bijin}, $\overline{\mathcal{F}(\beta)}\supset\overline{\mathcal{M}(\beta)}$. Combining this with Theorem \ref{th1} (1), we have $\overline{\mathcal{F}(\beta)}=\overline{\mathcal{M}(\beta)}$. Hence Theorem \ref{th1} (2) is obtained.

\vspace{0.2cm}
\noindent \textbf{Proof of Theorem \ref{th2} (1)}
\par (1) If $(\beta,\alpha)\notin\mathcal{M}$, it is clear that $I(\beta,\alpha)=\emptyset$. Let $(\beta,\alpha)\in\mathcal{M}$ and matching occurs at time $m-1$. Denote $k_{+}=(10a_{3}\cdots a_{m}\cdots)$ and $k_{-}=(01b_{3}\cdots b_{m}\cdots)$. By the proof of Theorem \ref{th1} (1), we have $(T^{+}_{\beta,\alpha})^{k}(c)=\beta^{k-2}\alpha+\beta^{k-3}\alpha-\beta^{k-3} a_{3}+\cdots +\beta\alpha-\beta a_{k-1}+\alpha-a_{k}$, for all $k\leq m$, and
$$\left \{ \begin{array}{ll}
(T^{+}_{\beta,\alpha})^{k}(c)<(1-\alpha)/\beta    \ \ \ \  \ \ {\rm if} \ a_{k+1}=0\\
(T^{+}_{\beta,\alpha})^{k}(c)\geq(1-\alpha)/\beta   \ \ \ \  \ \ {\rm if} \ a_{k+1}=1,
\end{array}
\right.
$$
where $(1-\alpha)/\beta$ represents the critical point $c$. Hence we can obtain $m$ inequalities of $\alpha$ in the following form:
$$\left \{ \begin{array}{ll}
\alpha<(1+\sum_{i=1}^{k-2}\beta^{i}a_{k-i+1})/(1+\beta+\cdots+\beta^{k-1})    \ \ \ \  \ \ {\rm if} \ a_{k+1}=0\\
\alpha\geq(1+\sum_{i=1}^{k-2}\beta^{i}a_{k-i+1})/(1+\beta+\cdots+\beta^{k-1})   \ \ \ \  \ \ {\rm if} \ a_{k+1}=1.
\end{array}
\right.
$$
Similarly, we have $ (T^{-}_{\beta,\alpha})^{k}(c)=\beta^{k-1}+\beta^{k-2}\alpha-\beta^{k-2}+\beta^{k-3}\alpha-\beta^{k-3}b_{3}+\cdots +\beta\alpha+\beta b_{k-1}+\alpha-b_{k},
$ and the following inequalities:

$$\left \{ \begin{array}{ll}
\alpha\leq(1+\sum_{i=3}^{k}b_{i}\beta^{k+1-i}+\beta^{k-1}-\beta^{k})/(1+\beta+\cdots+\beta^{k-1})    \ \ \ \  \ \ {\rm if} \ b_{k+1}=0\\
\alpha>(1+\sum_{i=3}^{k}b_{i}\beta^{k+1-i}+\beta^{k-1}-\beta^{k})/(1+\beta+\cdots+\beta^{k-1})   \ \ \ \  \ \ {\rm if} \ b_{k+1}=1.
\end{array}
\right.
$$
\par In this way, we obtain 2$m$ inequalities about $\alpha$, and denote $I^{\prime}(\beta,\alpha)$ as the intersection of all inequalities. By the construction of $I^{\prime}(\beta,\alpha)$, for each $\alpha^{\prime}\in I^{\prime}(\beta,\alpha)$, the kneading invariants of $(\beta,\alpha^{\prime})$ have the same first $m$ symbols with $(k_{+},k_{-})$. Moreover, in order to remain the $\beta$ unchanged, each $\alpha^{\prime}$ must have matching at time $m-1$. Hence we obtain $I(\beta,\alpha)=\{\beta\}\times I^{\prime}(\beta,\alpha)$. It can be seen that $I(\beta,\alpha)$ is a singleton or a subinterval of $\Delta(\beta)$.
$\hfill\square$

\vspace{0.2cm}

\par  By Proposition \ref{bijin}, for any $\epsilon>0$ and $(\beta,\alpha)\in I(\beta,\alpha)$, there exists $\alpha^{\prime}$ satisfying $|\alpha-\alpha^{\prime}|<\epsilon$, such that $(\beta,\alpha^{\prime})\in\mathcal{F}(\beta,\alpha)$. Hence $\mathcal{F}(\beta,\alpha)$ is dense in $I(\beta,\alpha)$ and Theorem \ref{th2} (2) is obtained.

\begin{figure}[htbp]
\centering
\includegraphics[width=0.8\textwidth]{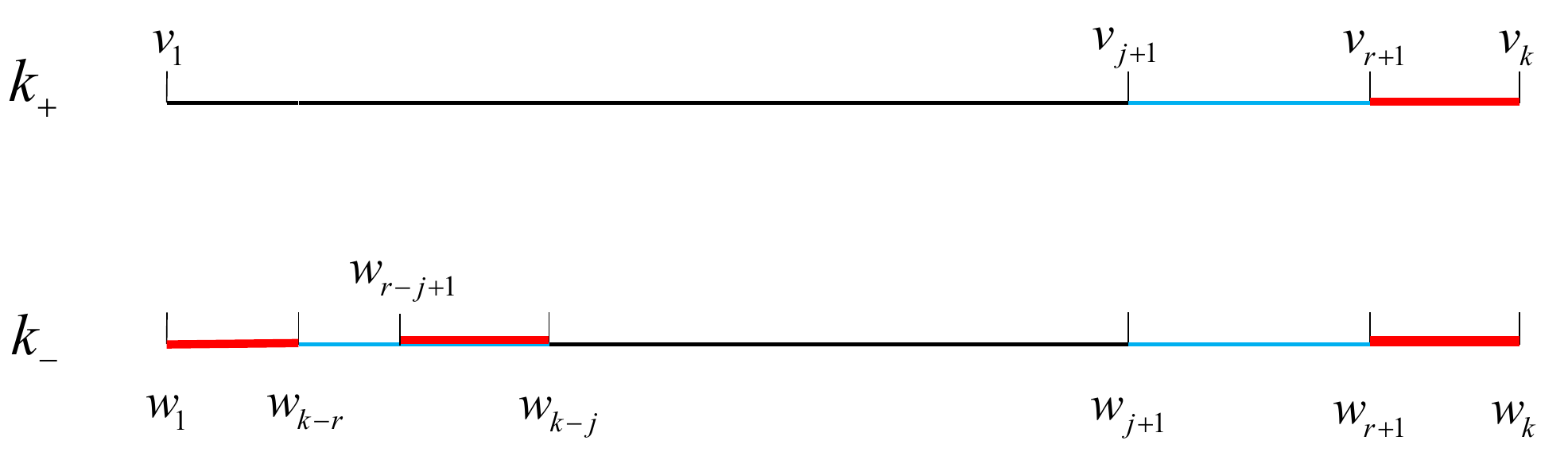}
\caption{Plot of Case A with $j<r<k-1$, four heavy lines represent the identical words.}
\label{fig:fig2}
\end{figure}

\begin{figure}[htbp]
\centering
\includegraphics[width=0.8\textwidth]{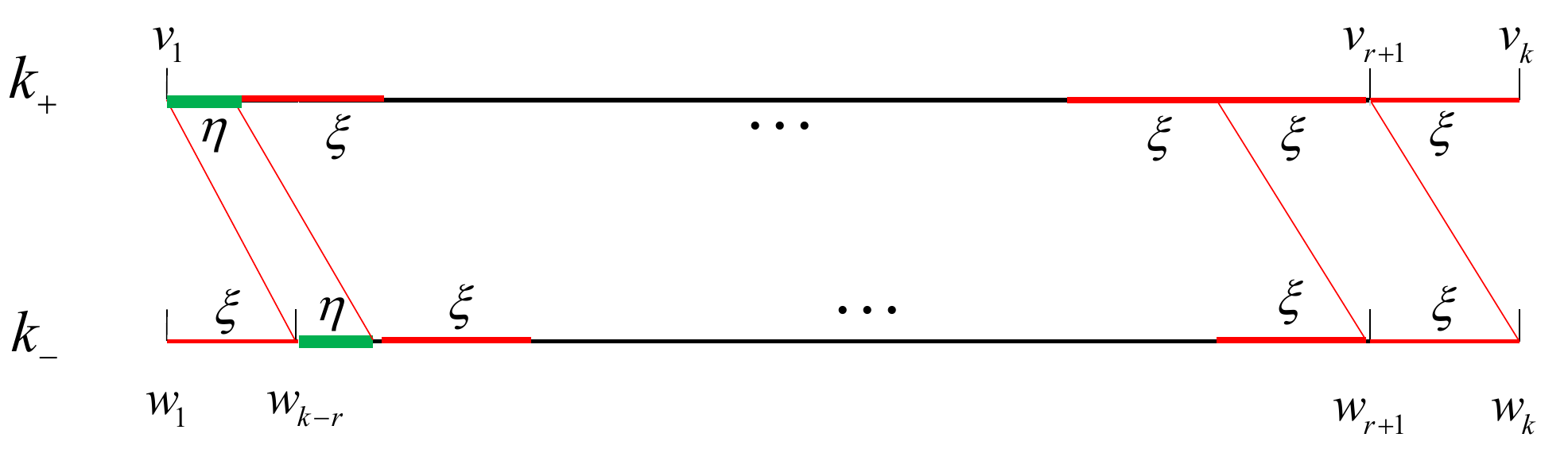}
\caption{Plot of Case B.}
\label{fig:fig2}
\end{figure}

\begin{proposition}\label{mononity} Let $(\beta,\alpha)$ and $(\beta,\alpha^{\prime})$ belong to the same fiber $\Delta(\beta)$, $(k_{+},k_{-})$ and $(k^{\prime}_{+},k^{\prime}_{-})$ be their corresponding  kneading invariants. Then $k^{\prime}_{+}\succ k_{+}\succ k^{\prime}_{-}\succ k_{-}$ if and only if $\alpha^{\prime}>\alpha$.

\begin{proof}
Denote $k_{+}=(10a_{3}\cdots a_{n}\cdots)$ and $k^{\prime}_{+}=(10a^{\prime}_{3}\cdots a^{\prime}_{n}\cdots)$, where $n$ be the maximal integer such that $a_{i}= a^{\prime}_{i}$, for all $i\leq n$. Since $k^{\prime}_{+}\succ k_{+}$, we have $a^{\prime}_{n+1}=1$ and $a_{n+1}=0$. By the proof of Theorem 1.3 (1), $a^{\prime}_{n+1}=1$ corresponds to $``\alpha^{\prime}\geq"$ and $a_{n+1}=0$ corresponds to $``\alpha<"$. Hence $\alpha^{\prime}>\alpha$ is obtained. Cases $k^{\prime}_{-}\succ k_{-}$, $k^{\prime}_{+}\preceq k_{+}$ and $k^{\prime}_{-}\preceq k_{-}$ can be obtained similarly.
\end{proof}
\end{proposition}

Before stating the following result, we give the definition of weak admissible. We say kneading invariants $(k_{+},k_{-})$ is weak admissible if satisfying
$$\sigma(k_{+}) \preceq\sigma^{n}(k_{+})\preceq\sigma(k_{-}) \ {\rm and} \ \sigma(k_{+})\preceq\sigma^{n}(k_{-})\preceq\sigma(k_{-}) \ \ \ \ {\rm for  \ all}\ \ n\geq 0,
$$
which means there may exist $n$ such that $\sigma^{n}(k_{+})=k_{-}$ or $\sigma^{n}(k_{-})=k_{+}$. Clearly, admissibility implies weak admissibility, but not vice versa. Now that we have known $I(\beta,\alpha)$
is a subinterval, then we only need to calculate the endpoints of $I(\beta,\alpha)$.
\begin{lemma}\label{calendpoint}
Let $T_{\beta,\alpha}$ has matching at time $m-1$, and the kneading invariants of $\Omega_{\beta,\alpha}$ be $(k_{+}=10\cdots a_{m}\cdots, k_{-}=01\cdots b_{m}\cdots)$. Then  left (right) endpoints of $I(\beta,\alpha)$ can be obtained by adding infinite sequence $\eta$ such that $(10\cdots a_{m}\eta, 01\cdots b_{m}\eta)$ is lexicographically smallest (biggest) weak admissible.
\begin{proof} By the definition of $I(\beta,\alpha)$, for all $(\beta,\alpha^{\prime})\in I(\beta,\alpha)$, their kneading invariants have the same first $m$ symbols. If there exists infinite sequence $\eta$ can be added after $a_{m}$ and $b_{m}$, such that $(10\cdots a_{m}\eta, 01\cdots b_{m}\eta)$ is weak admissible, and is meanwhile the lexicographically smallest with the fixed first $m$ symbols. Then by Proposition \ref{mononity}, the corresponding parameter pair $(\beta,\alpha_{l})$ is exactly the left endpoint. Similarly, if there exists infinite sequence $\xi$ can be added after $a_{m}$ and $b_{m}$, such that $(10\cdots a_{m}\xi, 01\cdots b_{m}\xi)$ is weak admissible, and is meanwhile lexicographically biggest, then the corresponding parameter pair $(\beta,\alpha_{r})$ is exactly the right endpoint. Moreover, if $\alpha_{l}$=$\alpha_{r}$, then $I(\beta,\alpha)$ is a singleton.
\par The key is to find out $\xi$ and $\eta$, here we only show the case of $\eta$. Denote $(k^{\prime}_{+},k^{\prime}_{-})=(10\cdots a_{m}\eta, 01\cdots b_{m}\eta)$. In order to obtain the lexicographically smallest weak admissible kneading invariants, we intend to add consecutive $0'$s after $a_{m}$ as much as possible. Since we want $(k^{\prime}_{+},k^{\prime}_{-})$ to be weak admissible, the number of consecutive $0'$s added after $a_{m}$ can not be greater than the number of consecutive $0'$s after $a_{1}=1$. If there is no $1$ among $a_{3}\cdots a_{m}$, then infinite $0'$s can be added and $k_{+}^{\prime}=10^{\infty}$. If there exists $1$ among $a_{3}\cdots a_{m}$, we discuss into two cases.
Case one, $a_{m}=0$ and $b_{m}=1$, we select the maximal integer $k$ ($k< m$) such that $a_{k+1}\cdots a_{m}=a_{1}\cdots a_{m-k}$, then it can be verified that $k_{+}^{\prime}=(10a_{3}\cdots a_{k})^{\infty}$. Case two, $a_{m}=1$ and $b_{m}=0$. At this case, select the maximal integer $r$ ($r<m$) such that $b_{r+1}\cdots b_{m}=a_{1}\cdots a_{r}$, then $k_{+}^{\prime}=10\cdots a_{m-r}(a_{m-r+1}\cdots a_{m})^{\infty}$. Since we only want to calculate the value of $\alpha_{l}$ or $\alpha_{r}$ , and use the formula $\alpha  = (\beta - 1)(K_{+}(1/\beta)-1)$ or $\alpha  = (\beta - 1)(K_{-}(1/\beta)-1)$, we do not need to ensure that $(k^{\prime}_{+},k^{\prime}_{-})$ is linearizable, just weak admissible is enough.  See Example \ref{endpoint} for an intuitive understanding. As a result, we can obtain the matching interval $I(\beta,\alpha)$ by only calculating the endpoints.
\end{proof}
\end{lemma}
\noindent \textbf{Proof of Theorem \ref{th2} (3)}
\par  When $(\beta,\alpha)\notin\mathcal{M}$, $I(\beta,\alpha)$ is empty and there is no endpoint. When $(\beta,\alpha)\in\mathcal{M}$ and $I(\beta,\alpha)$ is a singleton, we have $(\beta,\alpha_{l})=(\beta,\alpha_{r})=(\beta,\alpha)$ and by Remark \ref{tails}, $(\beta,\alpha)\in\mathcal{F}(\beta)$. When $I(\beta,\alpha)$ is a subinterval, $\overline{I(\beta,\alpha)}$ contains two different endpoints, the left endpoint and the right endpoint, denoted as $(\beta,\alpha_{l})$ and $(\beta,\alpha_{r})$. Here we only prove the case $(\beta,\alpha_{l})$, and the right endpoint can be proved similarly. Suppose $T_{\beta,\alpha}$ has matching at time $m-1$, and the kneading invariants of $\Omega_{\beta,\alpha}$ be $(k_{+}=10\cdots a_{m}\cdots, k_{-}=01\cdots b_{m}\cdots)$. By the proof of Lemma \ref{calendpoint}, we can obtain the lexicographically smallest weak admissible kneading invariants $(k^{\prime}_{+},k^{\prime}_{-})$. It is divided into 3 cases.
\par Case one, there is no $1$ among $a_{3}\cdots a_{m}$. We have $k^{\prime}_{+}=10^{\infty}$, and $k^{\prime}_{-}=01b_{3}\cdots b_{m-1}10^{\infty}$. It is clear that $(k^{\prime}_{+},k^{\prime}_{-})$ is weak admissible but not admissible, hence $(k_{+}^{\prime},k_{-}^{\prime})$ is not the real kneading invariants for $(\beta,\alpha_{l})$. Denote $k^{\prime\prime}_{-}=(01b_{3}\cdots b_{m-1})^{\infty}$, it can be verified that $K^{\prime}_{+}(t)=K^{\prime}_{-}(t)=K^{\prime\prime}_{+}(t)$, $(k^{\prime}_{+},k^{\prime\prime}_{-})$ is linearizable and corresponding to the left endpoint $(\beta,\alpha_{l})$. By Theorem \ref{lazysft}, at this case, $(\beta,\alpha_{l})\in\mathcal{F}(\beta)$.
\par Case two, there exists $1$ among $a_{3}\cdots a_{m}$, meanwhile $a_{m}=0$ and $b_{m}=1$. We can find the maximal integer $k<m$ such that $a_{k+1}\cdots a_{m}=a_{1}\cdots a_{m-k}$, $k_{+}^{\prime}=(10a_{3}\cdots a_{k})^{\infty}$ and meanwhile $k_{-}^{\prime}=01\cdots b_{m}(a_{m-k+1}\cdots $
$ a_{m})^{\infty}$. Denote $k_{-}^{\prime\prime}=(01\cdots b_{m}b_{m+1}$ $\cdots b_{2k-1})^{\infty}$, for the same reason,
$K^{\prime}_{+}(t)=K^{\prime}_{-}(t)=K^{\prime\prime}_{+}(t)$, $(k^{\prime}_{+},k^{\prime\prime}_{-})$ is linearizable and corresponding to the left endpoint $(\beta,\alpha_{l})$. At this case, both $k^{\prime}_{+}$ and $k^{\prime\prime}_{-}$ are periodic. Hence $(\beta,\alpha_{l})\in\mathcal{F}(\beta)$.
\par Case three, $a_{m}=1$ and $b_{m}=0$. Then we can find the maximal integer $r$ ($r<m$) such that $b_{r+1}\cdots b_{m}=a_{1}\cdots a_{r}$, $k_{+}^{\prime}=10\cdots a_{m-r}(a_{m-r+1}\cdots a_{m})^{\infty}$ and meanwhile $k_{-}^{\prime}=01\cdots b_{m-r}k_{+}^{\prime}$. Denote $k_{-}^{\prime\prime}=(01\cdots b_{r})^{\infty}$, similarly, $K^{\prime}_{+}(t)=K^{\prime}_{-}(t)=K^{\prime\prime}_{+}(t)$, $(k^{\prime}_{+},k^{\prime\prime}_{-})$ is linearizable and corresponding to the left endpoint $(\beta,\alpha_{l})$. Notice that $10\cdots a_{r}\neq a_{r+1}\cdots a_{m}$, otherwise matching will occur at time $r-1$, hence $k_{+}^{\prime}$ is eventually periodic and $(\beta,\alpha_{l})\in \mathcal{S}(\beta)\setminus\mathcal{M}(\beta)$.

 \par The result of right endpoint can also be obtained in the same way. If there is no $0$ among $b_{3}\cdots b_{m}$, we have $(\beta,\alpha_{r})\in \mathcal{F}(\beta)$. If there exists $0$ among $b_{3}\cdots b_{m}$ with $a_{m}=0$ and $b_{m}=1$, then $(\beta,\alpha_{r})\in \mathcal{F}(\beta)$.
If $a_{m}=1$ and $b_{m}=0$, then $(\beta,\alpha_{r})\in \mathcal{S}(\beta)\setminus\mathcal{M}(\beta)$. See Example \ref{endpoint} for intuitive explanation. Furthermore, the case of $I(\beta,\alpha)$ being a singleton is contained in Case one and Case two.
$\hfill\square$

\vspace{0.2cm}
\noindent \textbf{Proof of Corollary \ref{cor1.4}}
\par By the definition of $I(\beta,\alpha)$,
different matching intervals are disjoint. Based on Theorem \ref{th2} (3), we claim that even the closure of different matching intervals are also disjoint, that is, for any $I(\beta,\alpha_{1})\neq I(\beta,\alpha_{2})$, $\overline{I(\beta,\alpha_{1})}\cap \overline{I(\beta,\alpha_{2})}=\emptyset$. For convenience, we assume $I(\beta,\alpha_{1})$ lies on the left side of $I(\beta,\alpha_{2})$, denote $(\beta,\alpha_{r})$ as the right endpoint of $I(\beta,\alpha_{1})$ and $(\beta,\alpha_{l})$ as the left endpoint of $I(\beta,\alpha_{2})$.
 We prove this by contradiction, if $\overline{I(\beta,\alpha_{1})}\cap \overline{I(\beta,\alpha_{2})}\neq\emptyset$, then $(\beta,\alpha_{r})=(\beta,\alpha_{l})$. Denote the intersection as $(\beta,\alpha)$ and discuss into 3 cases.
\par Case one, $(\beta,\alpha)\notin I(\beta,\alpha_{1})$ and $(\beta,\alpha)\notin I(\beta,\alpha_{2})$. By the proof of Theorem \ref{th2} (3), both upper kneading sequence and lower kneading sequence of $(\beta,\alpha)$ are periodic, which indicates $(\beta,\alpha)\in\mathcal{F}$, and contradicts with the assumption that $(\beta,\alpha)$ not belonging to any of these two matching intervals.
\par Case two, $(\beta,\alpha)$ belongs to only one of the matching intervals. Such as, $(\beta,\alpha)\in I(\beta,\alpha_{1})$ but $(\beta,\alpha)\notin I(\beta,\alpha_{2})$. By Remark \ref{tails}, $(\beta,\alpha)\in I(\beta,\alpha_{1})$ indicates $(\beta,\alpha)\in \mathcal{F}$, $(\beta,\alpha)\notin I(\beta,\alpha_{1})$ implies $(\beta,\alpha)\in \mathcal{S}(\beta)\setminus\mathcal{M}(\beta)$, which leads to a contradiction.

\par Case three, $(\beta,\alpha)$ belongs to the both matching intervals. This indicates that any parameter pair in the two matching interval has the same matching, which contradicts with the assumption $I(\beta,\alpha_{1})\neq I(\beta,\alpha_{2})$.                           $\hfill\square$
\vspace{0.2cm}

\vspace{0.2cm}

\noindent \textbf{Proof of Corollary \ref{cor1.3}}
\par Let $\beta$ be a multinacci number of order $m$, then $\beta$ is the unique real solution to the equation $\beta^{m}-\beta^{m-1}-\cdots-1=0$ in the interval (1, 2). The necessity of this corollary was proved in \cite[Proposition 5.1]{bruin2017} and \cite[Proposition 1]{quackenbush2020periodic}, and can also be obtained by the proof of Theorem \ref{th2} (1). Observed that for any $\alpha\in(0,2-\beta)$, the kneading invariants of $\Omega_{\beta,\alpha}$ are
 \begin{equation}\label{m-times}
k_{+}|_{m+1}=1\underbrace{0\cdots0}_{m-times} \ \ \ {\rm and} \ \ \    k_{-}|_{m+1}=0\underbrace{1\cdots1}_{m-times}.
\end{equation}

\par Next we prove that if $\beta$ is not a multinacci number, then $I(\beta,\alpha)\neq\Delta(\beta)$ for any $\alpha\in(0,2-\beta)$. If $(\beta,\alpha)\notin\mathcal{M}$, $I(\beta,\alpha)=\emptyset$. If $T_{\beta,\alpha}$ has matching at time $m$, by Lemma \ref{zeros} and Proposition \ref{SFTIFF} , $\beta$ must be the biggest real root of the equation $a_{1}\beta^{m}-a_{2}\beta^{m-1}-\cdots-a_{m+1}=0$, where $a_{i}\in\{0,\pm1\}$ ($i\in\{1,\cdots,m+1\}$), $a_{1}=a_{2}=1$ and $a_{m+1}\neq0$. In addition, $\beta$ is not a multinacci number indicates that there exists at least one $i$ such that $a_{i}\neq1$. Hence the first $m+1$ symbols of $k_{+}$ and $k_{-}$ can not be totally symmetric as the equality (\ref{m-times}) above. However, by the proof of Lemma \ref{calendpoint},
$\alpha=0$ requires there exists no $1$ among $a_{2}\cdots a_{m+1}$, $\alpha=2-\beta$ requires there exists no $0$ among $b_{2}\cdots b_{m+1}$, which contradicts with the existence of $a_{i}\neq1$. Hence the endpoints of $I(\beta,\alpha)$ will not reach $0$ or $2-\beta$, then if $\beta$ is not multinacci number, $I(\beta,\alpha) \neq\Delta(\beta)$.
$\hfill\square$

\section{Some examples}
\begin{example}\label{example1}($(\beta,\alpha)\in\mathcal{M}\setminus\mathcal{F}$) \\
Let $k_{+}=100(10)^{\infty}$ and $k_{-}=011(10)^{\infty}$, it can be calculated that $\beta=(\sqrt{5}+1)/2$ and $\alpha\approx 0.23607$. By Proposition \ref{SFTIFF}, $T_{\beta,\alpha}$ has  matching at the second iteration, that is $T^{2}_{\beta,\alpha}(0)=T^{2}_{\beta,\alpha}(1)$. However, both $k_{+}$ and $k_{-}$ are not periodic implies $(\beta,\alpha)\notin\mathcal{F}$.

\end{example}

\begin{example}\label{realkneading}(Unique linearizable $(k_{+},k_{-})$ corresponds to $(\beta,\alpha)$)

\par (1) For the case $m=1$, let $k_{+}= (100011011)^{\infty}$ and $k_{-} = (011100)^{\infty}$, then $(k_{+},k_{-})$ can be non-periodic renormalized by $(w_{+}=100, w_{-}=011)$. Both $(k_{+},k_{-})$ and $(w_{+}^{\infty},w_{-}^{\infty})$ have the same parameter pair $(\beta,\alpha)\approx(1.618,0.191)$, but only $((100)^{\infty},(011)^{\infty})$ is the linearizable kneading invariants of the intermediate $\beta$-shift $\Omega_{\beta,\alpha}$.
\par (2) See the case $m>1$, let $k^{\prime}_{+}= (100101011010)^{\infty}$ and $k^{\prime}_{-} = (011010100101)^{\infty}$, it can be calculated that $(\beta,\alpha)\approx(1.272,0.364)$. $(k^{\prime}_{+},k^{\prime}_{-})$ can firstly be periodic renormalized via $(w_{+,1}=10,w_{-,1}=01)$, and then be non-periodic renormalized by $(w_{+,2}=100,w_{-,2}=011)$. By Lemma \ref{xunibeta}, only $(10,01)\ast((100)^{\infty},(011)^{\infty})=((100101)^{\infty},(011010)^{\infty})$ is the linearizable kneading invariants for parameter pair $(\beta,\alpha)$.
\end{example}

\begin{example}\label{manyk} (Cardinality of $\mathcal{K}^{\prime}$ is infinite)
\\
Let
$$
\left \{ \begin{array}{ll}
v=100011101101101101011(01)^{\infty},\\
w=011101101101101101011(01)^{\infty}.
\end{array}
\right.
$$
It can be seen that $(v,w)$ is linearizable and matching occurs at time 4. According to the selection of $k$ in the proof of Proposition \ref{bijin}, here we have $\mathcal{K}=\{10,13,16,18,21\}\cup\{21+2n, n\in\mathbb{N}\}$. Notice that when $k\in\{10,13,16\}$, $(k_{+},k_{-})$ is not linearizable, and for any $k\geq 18$, $(k_{+},k_{-})$ is  linearizable. Hence $\mathcal{K}^{\prime}=\mathcal{K}\setminus\{10,13,16\}$.
\end{example}

\begin{example}\label{emptybeta}($\mathcal{M}(\beta)=\emptyset$)\\
Let $k_{+}=100(01)^{\infty}$ and $ k_{-}=011(10)^{\infty}$, we can calculate that $\beta\approx1.76929$ and $\beta$ is a Perron number. By Lemma \ref{symmetric}, $\alpha=1-\beta/2$. However, the kneading determinant
$$K(t)=1-t-2t^2+2t^4=(t-1)(2t^3+2t^2-1),
$$
which means $2t^3+2t^2-1=0$ is the minimal polynomial of $t$. Notice that $t$ is not an algebraic integer, which contradicts with the fact that $1/\beta$ must be an algebraic integer if $(\beta,\alpha)\in\mathcal{M}$. Hence, in this case, $\mathcal{M}(\beta)=\emptyset$. However, we also have examples that $\beta$ is a Perron number but $\mathcal{M}(\beta)\neq\emptyset$, such as $k_{+}=(100010)^{\infty}$ and $k_{-}=(011101)^{\infty}$.
\end{example}

\begin{example}\label{unique}($I(\beta,\alpha)=\{(\beta,\alpha)\}$)\\
Let $k_{+}=(100)^{\infty}$ and $ k_{-}=(01)^{\infty}$, then $\beta\approx1.32472$, $\alpha\approx0.24512$. On the one hand, we have $I(\beta,\alpha)=\{(\beta,\alpha)\}$ via the proof of Theorem \ref{th2} (1). On the other hand, matching occurs at time $5$ and $a_{6}=0$, $b_{6}=1$.
By Lemma \ref{calendpoint}, both $0$ and $1$ can not be added, hence $I(\beta,\alpha)$ is a singleton.
\end{example}

\par There may exist many different matching intervals on $\Delta(\beta)$. By Lemma \ref{symmetric}, if there is a matching interval on the left of $1-\beta/2$, there must be a symmetric one on the other side.

\begin{example}\label{notunique}(Different $I(\beta,\alpha)$ on $\Delta(\beta)$)\\
Let $(k_{+}^{1},k_{-}^{1})=((1000)^{\infty},(01)^{\infty})$ and  $(k^{2}_{+},k^{2}_{-})=((100)^{\infty},(0110)^{\infty})$, then $(\beta,\alpha_{1})\approx(1.46557,$ $0.1288)$ and $(\beta,\alpha_{2})\approx(1.46557,0.2168)$, respectively. Their matching intervals are different since matching time is not identical. By  Theorem \ref{th2} (1), we have $I(\beta,\alpha_{1})=\{\beta\}\times(0, 0.1288]$ and $I(\beta,\alpha_{2})=\{(\beta,\alpha_{2})\}$. According to Lemma \ref{symmetric}, there also exist two distinct matching intervals induced by $((1001)^{\infty},(011)^{\infty})$ and $((10)^{\infty},(0111)^{\infty})$, denoted as $I(\beta,\alpha_{3})=\{(\beta,\alpha_{3}\approx0.3177)\}$ and $I(\beta,\alpha_{4})=\{\beta\}\times[0.4056,0.5344)$, respectively.
\end{example}

\begin{example}\label{endpoint}(Endpoints of $I(\beta,\alpha)$)\\
Here we give three examples to correspond with three cases in the proof of Theorem \ref{th2} (3).
\par (1) Let $k_{+}=(1000)^{\infty}$ and $k_{-}=(01)^{\infty}$, then matching occurs at time $3$ and there is no $1$ in $a_{2}a_{3}a_{4}=000$. By Theorem \ref{th2} (3), the kneading invariants of $\Omega_{\beta,\alpha_{l}}$ is $(10^{\infty},(010)^{\infty})$. Hence $\alpha_{l}=0$ and $(\beta,\alpha_{l})\in\mathcal{F}$.

\par (2) Let $k_{+}=(100010)^{\infty}$ and $k_{-}=(011111)^{\infty}$, then matching occurs at time $5$ and $a_{6}=0$, $b_{6}=1$. By Theorem \ref{th2} (3), the kneading invariants of $\Omega_{\beta,\alpha_{l}}$ is $((1000)^{\infty},(01111100)^{\infty})$. Hence $(\beta,\alpha_{l})\in\mathcal{F}$, and meanwhile $(\beta,\alpha_{l})\in I(\beta,\alpha)$.

\par (3) Let $k_{+}=(10001)^{\infty}$ and $k_{-}=(01110)^{\infty}$, then matching occurs at time $4$ and $a_{5}=1$, $b_{5}=0$. By Theorem \ref{th2} (3), the kneading invariants of $\Omega_{\beta,\alpha_{l}}$ is $(10(001)^{\infty},(011)^{\infty})$. Hence $(\beta,\alpha_{l})\in\mathcal{S}\setminus\mathcal{M}$, but $(\beta,\alpha_{l})\notin I(\beta,\alpha)$.
\end{example}

\section{Final comments}
The results of \cite[Theorem 3]{glendinning1996zeros} and \cite[Theorem 1.3]{li2016intermediate} show that if $\beta$ is a transcendental number, then $\mathcal{F}(\beta)=\emptyset$, by Remark 1.2 (2), hence $\mathcal{M}(\beta)=\emptyset$.  Indeed, for $\Omega_{\beta, \alpha}$ to be a SFT, we require $\beta \in (1, 2)$ to be a maximal root of a polynomial with coefficients in $\{ -1, 0, 1 \}$.  Moreover, the entropy of a SFT is the logarithm of the largest eigenvalue $\lambda$ of a nonnegative integral matrix.  By \cite[Theorem 3]{lind1984entropies}, we have that $\lambda$ is the positive $n^{\textup{th}}$-root of a Perron number.  Since the entropy of an intermediate $\beta$-shift is $\ln(\beta)$, if $\mathcal{F}(\beta)\neq\emptyset$, then there exists an $n \in \mathbb{N}$ such that $\beta$ is the positive $n^{\textup{th}}$-root of a Perron number. However, there exist lots of examples such that $\mathcal{M}(\beta)=\emptyset$ when $\beta$ is a Perron number, see Example \ref{emptybeta}. A natural question is that, can we classify the Perron numbers according to $\mathcal{M}(\beta)=\emptyset$ or $\mathcal{M}(\beta)\neq\emptyset$? This is an interesting question still unknown.

\par On the other side, when considering the relationships between sofic, SFT and matching, we have $\overline{\mathcal{F}(\beta)}=\overline{\mathcal{M}(\beta)}=\overline{\mathcal{S}(\beta)}=\overline{\Delta(\beta)}$ in the case $\beta$ being a multinacci number. It seems that the relationships also hold for $\beta$ being other Pisot number, but still unknown. And the key is to prove that $\mathcal{M}(\beta)$ is dense in $\Delta(\beta)$.

\bibliographystyle{plain}
\bibliography{subshifts-of-finite-type-and-Matching}

\begin{thebibliography}{10}

\bibitem{alseda1996}
Ll~Alsed{\`a} and F~Manosas.
\newblock Kneading theory for a family of circle maps with one discontinuity.
\newblock {\em Acta Math. Univ. Comenian.(NS)}, 65(1):11--22, 1996.

\bibitem{barnsley2012}
Michael Barnsley, Brendan Harding, and Andrew Vince.
\newblock The entropy of a special overlapping dynamical system.
\newblock {\em Ergodic Theory and Dynamical Systems}, 34(2):483--500, 2012.

\bibitem{barnsley2014}
Michael Barnsley, Wolfgang Steiner, and Andrew Vince.
\newblock Critical itineraries of maps with constant slope and one
  discontinuity.
\newblock In {\em Mathematical Proceedings of the Cambridge Philosophical
  Society}, volume 157, pages 547--565. Cambridge University Press, 2014.

\bibitem{blanchard1989}
Fran{\c{c}}ois Blanchard.
\newblock $\beta$-expansions and symbolic dynamics.
\newblock {\em Theoretical Computer Science}, 65(2):131--141, 1989.

\bibitem{bruin2017}
Henk Bruin, Carlo Carminati, and Charlene Kalle.
\newblock Matching for generalised $\beta$-transformations.
\newblock {\em Indagationes Mathematicae}, 28(1):55--73, 2017.

\bibitem{cui2015}
Hongfei Cui and Yiming Ding.
\newblock Renormalization and conjugacy of piecewise linear lorenz maps.
\newblock {\em Advances in Mathematics}, 271:235--272, 2015.

\bibitem{ding2011}
Yiming Ding.
\newblock Renormalization and $\alpha$-limit set for expanding lorenz maps.
\newblock {\em Discrete \& Continuous Dynamical Systems}, 29(3):979, 2011.

\bibitem{ding2021complete}
Yiming Ding and Yun Sun.
\newblock Complete invariants and parametrization of expansive lorenz maps.
\newblock {\em arXiv preprint arXiv:2103.16979}, 2021.

\bibitem{ding2022alpha}
Yiming Ding and Yun Sun.
\newblock $\alpha$-limit sets and lyapunov function for maps with one
  topological attractor.
\newblock {\em Acta Mathematica Scientia}, 42(2):813--824, 2022.

\bibitem{glendinning1990topological}
Paul Glendinning.
\newblock Topological conjugation of lorenz maps by $\beta$-transformations.
\newblock In {\em Mathematical Proceedings of the Cambridge Philosophical
  Society}, volume 107, pages 401--413. Cambridge University Press, 1990.

\bibitem{glendinning1996zeros}
Paul Glendinning and Toby Hall.
\newblock Zeros of the kneading invariant and topological entropy for lorenz
  maps.
\newblock {\em Nonlinearity}, 9(4):999, 1996.

\bibitem{glendinning1993prime}
Paul Glendinning and Colin Sparrow.
\newblock Prime and renormalisable kneading invariants and the dynamics of
  expanding lorenz maps.
\newblock {\em Physica D: Nonlinear Phenomena}, 62(1-4):22--50, 1993.

\bibitem{hubbard1990classification}
John~H Hubbard and Colin~T Sparrow.
\newblock The classification of topologically expansive lorenz maps.
\newblock {\em Communications on Pure and Applied Mathematics}, 43(4):431--443,
  1990.

\bibitem{kalle2012beta}
Charlene Kalle and Wolfgang Steiner.
\newblock Beta-expansions, natural extensions and multiple tilings associated
  with pisot units.
\newblock {\em Transactions of the American Mathematical Society},
  364(5):2281--2318, 2012.

\bibitem{komornik2011}
Vilmos Komornik and P~Loreti.
\newblock Expansions in noninteger bases.
\newblock {\em Integers}, 11(A9):30, 2011.

\bibitem{komornik1998}
Vilmos Komornik and Paola Loreti.
\newblock Unique developments in non-integer bases.
\newblock {\em The American mathematical monthly}, 105(7):636--639, 1998.

\bibitem{kraaikamp2012natural}
Cor Kraaikamp, Thomas~A Schmidt, and Wolfgang Steiner.
\newblock Natural extensions and entropy of $\alpha$-continued fractions.
\newblock {\em Nonlinearity}, 25(8):2207, 2012.

\bibitem{li2016intermediate}
Bing Li, Tuomas Sahlsten, and Tony Samuel.
\newblock Intermediate $\beta$-shifts of finite type.
\newblock {\em Discrete \& Continuous Dynamical Systems}, 36(1):323--344, 2016.

\bibitem{li2019denseness}
Bing Li, Tuomas Sahlsten, Tony Samuel, and Wolfgang Steiner.
\newblock Denseness of intermediate $\beta$-shifts of finite-type.
\newblock {\em Proceedings of the American Mathematical Society},
  147(5):2045--2055, 2019.

\bibitem{lind1984entropies}
Douglas~A Lind.
\newblock The entropies of topological markov shifts and a related class of
  algebraic integers.
\newblock {\em Ergodic Theory and Dynamical Systems}, 4(2):283--300, 1984.

\bibitem{milnor1988iterated}
John Milnor and William Thurston.
\newblock On iterated maps of the interval.
\newblock In {\em Dynamical systems}, pages 465--563. Springer, 1988.

\bibitem{palmer1979classification}
Marion~R Palmer.
\newblock {\em On the classification of measure preserving transformations of
  Lebesgue spaces}.
\newblock PhD thesis, University of Warwick, 1979.

\bibitem{parry1960beta}
William Parry.
\newblock On the $\beta$-expansions of real numbers.
\newblock {\em Acta Mathematica Hungarica}, 11(3-4):401--416, 1960.

\bibitem{parry1964representations}
William Parry.
\newblock Representations for real numbers.
\newblock {\em Acta Mathematica Academiae Scientiarum Hungarica},
  15(1):95--105, 1964.

\bibitem{quackenbush2020periodic}
Blaine Quackenbush, Tony Samuel, and Matt West.
\newblock Periodic intermediate $\beta$-expansions of pisot numbers.
\newblock {\em Mathematics}, 8(6):903, 2020.

\bibitem{renyi1957representations}
Alfr{\'e}d R{\'e}nyi.
\newblock Representations for real numbers and their ergodic properties.
\newblock {\em Acta Math. Acad. Sci. Hungar}, 8(3-4):477--493, 1957.

\bibitem{sidorov2003}
Nikita Sidorov.
\newblock Almost every number has a continuum of $\beta$-expansions.
\newblock {\em The American Mathematical Monthly}, 110(9):838--842, 2003.

\end{thebibliography}

\end{document}